\documentclass{article}

\usepackage[preprint]{neurips_2020}
\usepackage[utf8]{inputenc} % allow utf-8 input
\usepackage[T1]{fontenc}    % use 8-bit T1 fonts
\usepackage{hyperref}       % hyperlinks
\usepackage{url}            % simple URL typesetting
\usepackage{booktabs}       % professional-quality tables
\usepackage{amsfonts}       % blackboard math symbols
\usepackage{nicefrac}       % compact symbols for 1/2, etc.
\usepackage{microtype}      % microtypography
\usepackage{amsmath}
\usepackage{bigints}
\usepackage{amsthm}
\usepackage{xpatch}
\usepackage{graphicx}
\usepackage{wrapfig}
\usepackage{verbatim}
\usepackage{multicol}
\usepackage{xcolor}
\makeatletter
   \xpatchcmd{\@thm}{\fontseries\mddefault\upshape}{}{}{} % sam
\newtheorem{prop}{Proposition}
\numberwithin{prop}{subsection}
\newtheorem{coro}{Corollary}
\numberwithin{coro}{subsection}
\newtheorem{defn}{Definition}
\numberwithin{defn}{subsection}
\usepackage{appendix}
\usepackage{bm}

\title{STENCIL-NET: Data-driven solution-adaptive discretization of partial differential equations}

\author{%
  Suryanarayana Maddu$^{1,2,3}$, Dominik Sturm$^{1,2,3}$,
  Bevan L. Cheeseman$^{4}$, \And  Christian L.~M\"{u}ller$^{5}$, Ivo F.~Sbalzarini$^{1,2,3}$ 
  \thanks{\scriptsize $^1$  Institute of Artificial Intelligence, Faculty of Computer Science, Technische Universi\"{a}t Dresden, Dresden, Germany. \scriptsize $^2$ Max Planck Institute of Molecular Cell Biology and Genetics, Center for Systems Biology Dresden, Dresden, Germany. \scriptsize $^3$ Cluster of Excellence Physics of Life, TU Dresden, Germany. \scriptsize $^4$ ONI, Oxford, United Kingdom. \scriptsize $^5$ Flatiron Institute, New York City, USA.}\\  \\  
}

\begin{document}
\maketitle
\begin{abstract}
Numerical methods for approximately solving partial differential equations (PDE) are at the core of scientific computing. Often, this requires high-resolution or adaptive discretization grids to capture relevant spatio-temporal features in the PDE solution, e.g., in applications like turbulence, combustion, and shock propagation. Numerical approximation also requires knowing the PDE in order to construct problem-specific discretizations. Systematically deriving such solution-adaptive discrete operators, however, is a current challenge. Here we present STENCIL-NET, an artificial neural network architecture for data-driven learning of problem- and resolution-specific local discretizations of nonlinear PDEs. STENCIL-NET achieves numerically stable discretization of the operators in an unknown nonlinear PDE by spatially and temporally adaptive parametric pooling on regular Cartesian grids, and by incorporating knowledge about discrete time integration. Knowing the actual PDE is not necessary, as solution data is sufficient to train the network to learn the discrete operators. A once-trained STENCIL-NET model can be used to predict solutions of the PDE on larger spatial domains and for longer times than it was trained for, hence addressing the problem of PDE-constrained extrapolation from data. To support this claim, we present numerical experiments on long-term forecasting of chaotic PDE solutions on coarse spatio-temporal grids. We also quantify the speed-up achieved by substituting base-line numerical methods with equation-free STENCIL-NET predictions on coarser grids with little compromise on accuracy.
\end{abstract}

\section{Introduction}

\begin{figure}
    \centering
    \includegraphics[width=1.0\textwidth]{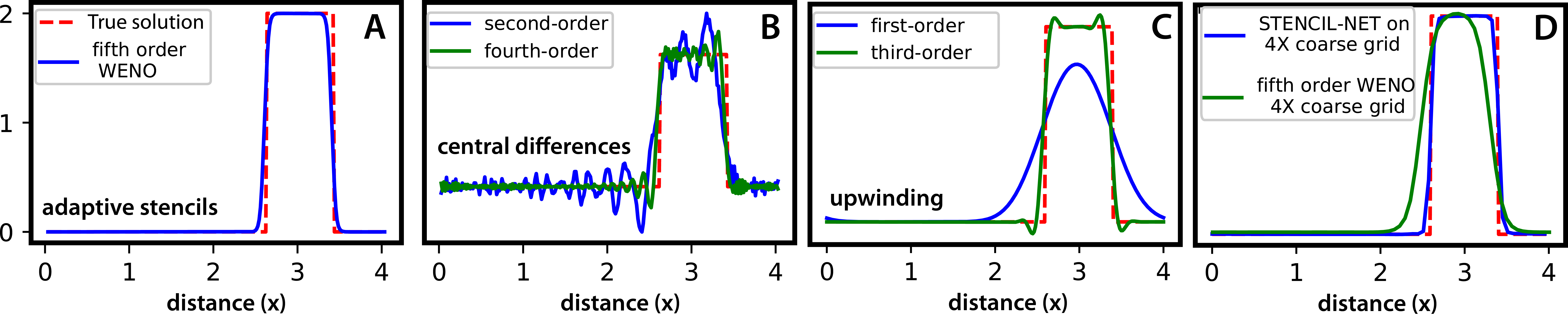}
     \vspace{-1.5em}
    \caption{ \footnotesize Numerical solutions of advection of a sharp pulse using different discretization schemes. The advection velocity is $c=2$ in the advection equation $u_t + cu_x = 0$. All solutions are visualized at time $t=3$\,s. 
    }
    \vspace{-1.5em}
    \label{steep}
\end{figure}

Mathematical modeling of spatio-temporal systems as Partial Differential Equations (PDEs) has been invaluable in gaining mechanistic insight from real-world dynamical processes. This has been further fueled by rapid advancements in numerical methods and parallel computing, enabling the study of increasingly complex systems in real-world geometries. Numerical methods like finite-difference (FD), finite-volume (FV), and finite-element (FE) methods are routinely used to approximate the solutions of PDEs on computational grids or meshes (\cite{patankar2018numerical,moin2010fundamentals}). Such mesh/grid-based methods rely on interpolating (in some basis) the PDE solution on discrete points (\cite{moin2010fundamentals}). This often requires high-resolution or adaptive-resolution grids with problem-specific discrete operators in order to accurately and stably represent small scales and high frequencies in the PDE solution, as illustrated in Fig.~\ref{steep}. 

Popular numerical schemes, such as Total Variation Diminishing (TVD), Essentially Non-Oscillatory (ENO), and Weighted Essentially Non-Oscillatory (WENO) methods use solution-adaptive nonlinear stencil weights in order to provide stable numerical solutions  (\cite{osher2004level,shu1998essentially}). Increasing the accuracy in such schemes, however, requires drastically higher-resolution spatial grids, incurring a large computational cost especially in three-dimensional simulations. 

Indeed, the computational cost for numerically solving advection-dominated problems (like turbulence, combustion, or shock propagation) using WENO schemes scales at least with $\Delta x^{-4}$, where $\Delta x$ is the spatial grid resolution. For this reason, and despite advancements in numerical algorithms and computer hardware, simulations of PDEs are yet to reach the speeds required for real-time applications. This is particularly important for optimization analysis applications, like uncertainty quantification and Bayesian inverse inference, which require repeated computation of forward solutions of PDEs (\cite{bijl2013uncertainty,mishra2018machine}) for different coefficients, benefiting from solution-specific discretization schemes at runtime. Therefore, there has been much recent interest in using modern machine-learning methods as predictive surrogate models for fast and accurate discretized PDE solution. They achieve speed-up by substituting compute-intensive parts of numerical solvers (\cite{tompson2017accelerating}) or by projecting the data to reduced representations using autoencoders for latent time integration (\cite{kim2019deep}). Alternatively, model-free methods based on Recurrent Neural Networks (RNNs) (\cite{pathak2017using,patankar2018numerical,vlachas2018data}) have shown unprecedented success in forecasting chaotic spatio-temporal dynamics. All these approaches, however, depend on knowing the actual PDE beforehand or on projection to a latent representation tailored to a specific domain size and geometry of the problem.

Here, we present an approach that does not require knowing the PDE and is purely data-driven. We present the STENCIL-NET architecture for data-driven,  solution-adaptive discretization of unknown nonlinear PDEs. STENCIL-NET uses a ``Network In Network'' (NIN) (\cite{lin2013network}) architecture to learn high-level discretization features from local input patches with reduced data requirements. The NIN network works by sliding a Multi-Layer Perceptron (MLP) over the input patches to perform cascaded cross-channel parametric pooling that enables learning complex solution features. We also draw formal connections between \textit{WENO-like} solution-adaptive numerical schemes, rational functions, and MLPs in order to theoretically defend our architecture choice.

\section{Background and related work}
Recently, there has been a surge of interest in using machine learning for discovery, prediction, and control of dynamical systems. Methods based on sparse-regression techniques have been proposed to learn symbolic ODEs/PDEs directly from limited amount of noise-free synthetic (\cite{brunton2016discovering,rudy2017data}), or noisy experimental data (\cite{maddu2019stability}). In the data-rich regime, \textit{Physics Informed Neural Networks} (PINN) have gained popularity for both discovery and approximate solution of PDEs (\cite{raissi2017physics}). PINNs use knowledge from physics in order to provide structural regularization to deep neural networks, and they use automatic differentiation to reconstruct the differential operators in the ODE/PDE models. Such physics-informed machine learning models have also been extended to reconstructing hidden latent variables (\cite{raissi2020hidden}) and to separating noise from signal using time-stepping constraints (\cite{raissi2018multistep,rudy2019deep}). Similarly, Convolutional Neural Networks (CNN) have been used for model discovery by learning discrete differential operators from data via sum-of-filter constraints on the filter weights (\cite{long2017pde,long2019pde}). 

Solution-adaptive discretizations of known nonlinear PDEs have been achieved using MLPs that learn problem-specific optimal parameters for classical numerical schemes on coarser grids (\cite{mishra2018machine}). The resulting consistent discretization on coarser grids can then be used to accelerate numerical simulations. 
Recent work on solution- and resolution-specific discretization of nonlinear PDEs has shown CNN filters that are able to generalize to larger spatial solution domains than what they have been trained on (\cite{bar2019learning}). However, an over-complete set of CNN filters was used, which renders learning high-level discretization features computationally expensive and data-demanding. In all previous works, the underlying PDE and the associated discrete numerical fluxes for learning accurate solution-specific discretizations had to be known beforehand. 

\section{Differential operators and convolutions} \label{Section3}
PDEs are among the most common mathematical models for  spatio-temporal dynamical processes. They have found extensive applications ranging from modeling ocean currents (\cite{higdon2006numerical}), to combustion kinetics (\cite{veynante2002turbulent}), to building predictive models of cell and tissue mechanics during growth and development (\cite{prost2015active}). In general, the spatio-temporal evolution of a state variable $u \in \mathbb{R}^d$ is given by a PDE of the form
\begin{equation}\label{PDE}
    \alpha_1 \frac{\partial u (x,t)}{\partial t}  + \alpha_2 \frac{\partial^2 u (x,t)}{\partial t^2} = \mathcal{N}\left( \Xi, u, \partial_x u, \partial_x f(u), \partial_{xx} u, \partial_{xxx} u, u^2, \ldots  \right), 
\end{equation}
where $\mathcal{N}\left( \cdot \right)$ is the nonlinear function that models the dynamics of a process, $\Xi$ is the set of PDE coefficients (e.g., diffusion constants or viscosities), and $f(u)$ are the flux terms, e.g., $f(u) = u^2, f(u)= cu$. For the time derivative, we only consider the binary and mutually exclusive case, i.e., $(\alpha_1 = 0, \alpha_2 = 1)$ or $(\alpha_1 = 1, \alpha_2 = 0)$.

There exists no one general method for numerically solving nonlinear PDEs. The most common methods like FD, FV, and FE work by discretizing the spatio-temporal domain by computational nodes called grid or mesh points, followed by solving the so-discretized right-hand side $\mathcal{N}_d (\cdot)$ with appropriate initial and/or boundary conditions. Here, we are only concerned with symmetric stencil configurations on uniformly spaced Cartesian grids. For clarity of presentation, we also derive our ideas in 1D, although they also work in higher dimensions. 

\subsection{Fixed-stencil approximations of differential operators}
The $l^{th}$ spatial derivative at location $x_i$ on a grid with spacing $\Delta x$ can be approximated with convergence order $r$ using linear convolutions:
\begin{equation}\label{linear_derv}
    \left.\frac{\partial^{l} u}{\partial x^{l}}\right \vert_{x=x_i} = \sum_{x_j \in S_m(x_i)} \!\!\!
    \xi_j u_j + \mathcal{O}(\Delta x^r), \:\:\:
\end{equation}
where $u_j = u(x_j)$. The $\xi_j$ are the stencil weights that can be determined by local polynomial interpolation  (\cite{fornberg1988generation,shu1998essentially}). The  stencil is $\mathcal{S}_m(x_i) = \{ x_{i-m}, x_{i-m+1},\ldots ,x_i, \ldots , x_{i+m-1}, x_{i+m} \}$ with size $\vert \mathcal{S}_m (x_i) \vert = 2m+1$. For a spatial domain of size $L$, the spacing is $\Delta x = L/N_{x}$, where $N_x$ is the number of grid points discretizing space. The following proposition from \cite{schrader2010discretization,long2017pde} defines the discrete moment conditions that need be fulfilled in order to obtain an $\mathcal{O}(\Delta x^r)$ accurate differentiation scheme.

\begin{prop} [\textbf{Discrete moment conditions}]  \label{prop1}
For an $r^{th}$ order approximation of the $l^{th}$ derivative at point $x_i$, the stencil needs to satisfy the following moment conditions:
\begin{equation*} 
    \sum_{x_j \in S_m(x_i)} \xi_j \left(\Delta x_{j}\right)^{k} = \begin{cases} 0 & 0 \leq k \leq l + r-1, k \neq l \\
(-1)^{k}k! & k = l \\
< \infty &  \textrm{otherwise,} \end{cases},
\end{equation*}
where $k,l,r \in \mathbb{Z}_{0}^{+}$ and $\Delta x_j = (x_i - x_j)$. A stencil fulfilling these conditions must have size $ \vert S_m \vert \geq \vert l+r-1 \vert $.
\end{prop}

\begin{prop} [\textbf{Nonlinear convolutions}] \label{prop2} Nonlinear terms of order 2, including product of derivatives (e.g. $u u_x , u^2, u_x u_{xx}$) can be approximated by nonlinear convolution. For $l_1,l_2, r_1, r_2 \in \mathbb{Z}_{0}^{+}$, we can write,
\begin{equation} \label{non_linear_convolutions}
   \left. \left(\frac{\partial^{l_1} u}{\partial x^{l_1}} \right)  \left(\frac{\partial^{l_2} u}{\partial x^{l_2}} \right) \right\vert_{x=x_i}  = \sum_{x_{j} \in S_{m} } \sum_{x_{k} \in S_{m}} \xi_{jk} u_j u_k + \mathcal{O}(\Delta x^{r_1 + r_2}), %\approx \textbf{p}(\mathbf{u}_m(x_i)),
\end{equation}
such that the stencil size $\vert S_m \vert \geq \vert l_1 + r_1 - 1\vert$, and $\vert S_m \vert \geq \vert l_2 + r_2 - 1\vert$ following proposition (\ref{prop1}). Eq. (\ref{non_linear_convolutions}) is a convolution with a  Volterra-based quadratic form, as described by \cite{zoumpourlis2017non}.
\end{prop}

\begin{wrapfigure}{r}{0.5\textwidth} 
\vspace{-2.0em}
\centering
\includegraphics[width=0.45\textwidth]{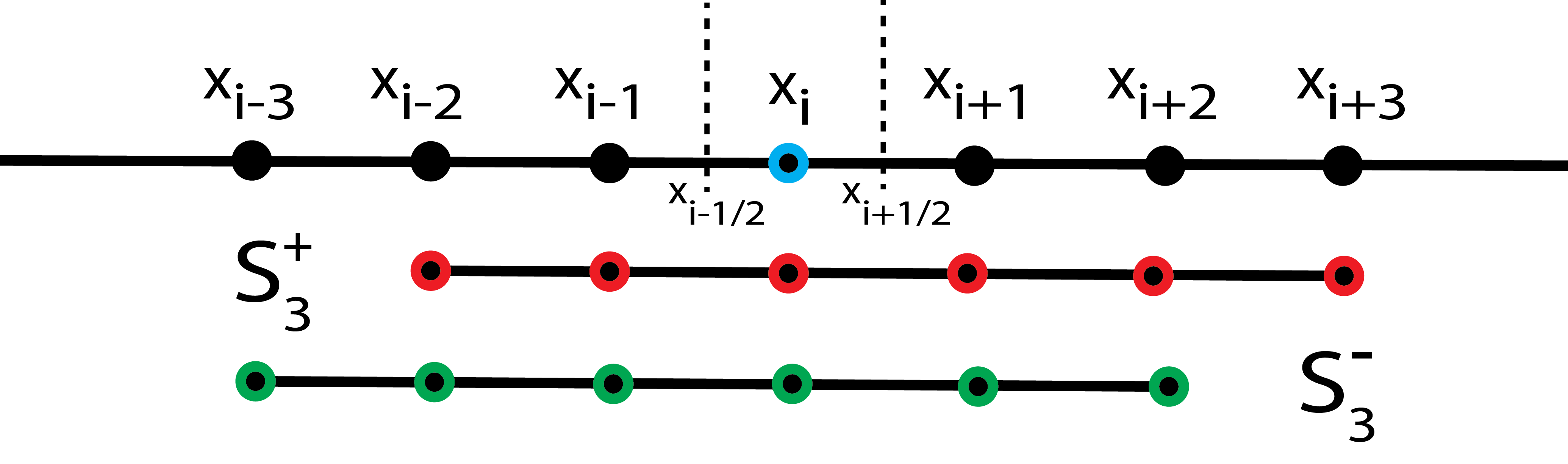}
\vspace{-0.3em}
\caption{\footnotesize Schematic of the stencil $S_3(x_i)$ centered around the point $x_i$. Here, $S_{3}(x_i) = S_3^{+} \cup S_3^{-}$.}
\label{stencil_draw}
\vspace{-1.0em}
\end{wrapfigure}

The linear convolutions in Eq.~(\ref{linear_derv}) and the nonlinear convolutions in Eq.~(\ref{non_linear_convolutions}) are based on local polynomial interpolation. They are thus useful to discretize smooth solutions arising in problems like reaction-diffusion, fluid flow at low Reynolds numbers, and elliptic PDEs. But they fail to discretize nonlinear flux terms arising, e.g., in problems like advection-dominated flows, Euler equations, multi-phase flows, level-set, and Hamilton Jacobi equations (\cite{osher2004level}). This is because higher-order interpolation near discontinuities leads to oscillations that do not decay when refining the grid, a fact known as ``Gibbs phenomenon''. Moreover, fixed-stencil weights computed from moment conditions can fail to capture the direction of information flow (``upwinding''), leading to non-causal and hence unstable discretizations (\cite{courant1952solution,patankar2018numerical}). The latter can be relaxed by biasing the stencil along the direction of information flow or by constructing numerical fluxes with smooth flux-splitting methods like Godunov and Lax-Friedrichs schemes. To counter the issue of spurious Gibbs oscillations, artificial viscosity (\cite{sod1985numerical}) can be added at the expense of solution accuracy. Alternatively, solution-adaptive stencils with ENO (Essentially Non-Oscillatory) and WENO (Weighted ENO) weights are available (\cite{shu1998essentially}).

\subsection{Motivation for solution-adaptive stencils: ENO/WENO approximations}
The ENO/WENO method for discretely approximating a continuous function $f$ at a point $x_{i \pm 1/2}$ can be written as a linear convolution:
\begin{equation}\label{WENO_form}
    \hat{f}_{i\pm {1/2}} = \sum_{x_j \in S^{\pm}(x_i)} \nu_j^{\pm} f(x_j) + \mathcal{O}(\Delta x^r),
\end{equation}
where $f_{i\pm 1/2} = f(x_i \pm \Delta x/2)$, and the function value and coefficient on the stencils, given by $\mathbf{f}_m = \{f(x_j): x_j \in S_m(x_i)\}$ and $\nu(x_j)$, respectively. The stencils $S_m^{\pm} (x_i)$, as illustrated in Fig.~\ref{stencil_draw}, are $S_m^{\pm}(x_i) = S_m(x_i) \setminus x_{i\pm m} $. However, unlike the fixed-stencil convolutions in Eqs.~(\ref{linear_derv}) and (\ref{non_linear_convolutions}), the coefficients ($\nu$) in ENO/WENO approximations are computed based on local smoothness indicators on smaller sub-stencils, which in turn depend on the functions values, leading to locally solution-adaptive weights (Proposition \ref{prop3}). This allows accurate function reconstruction even if the function $f$ is not smooth. 

\begin{figure}
    \centering
    \includegraphics[width=0.8\textwidth]{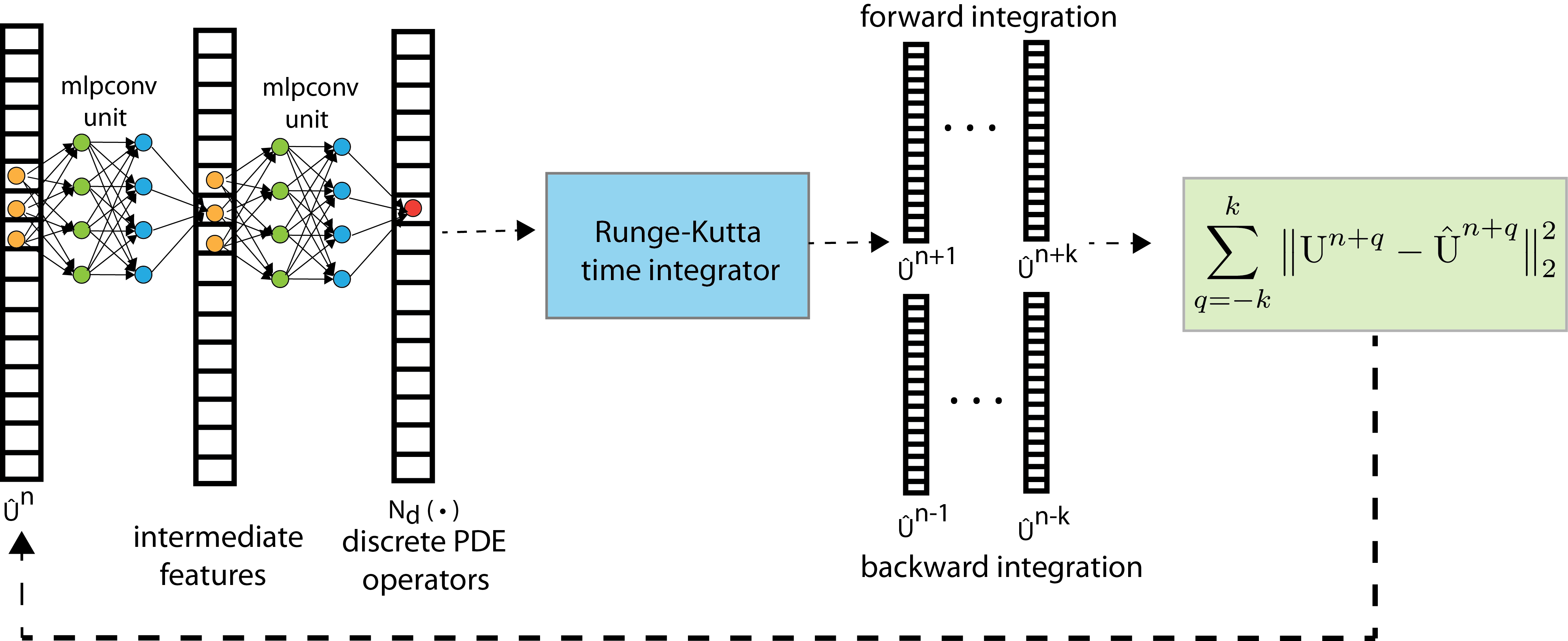}
    \caption{\footnotesize \textbf{STENCIL-NET for data-driven discretization:} The mlpconv unit performs parametric pooling by moving the MLP network across the input vector $\hat{\mathbf{U}}^n$ at time $n$ to generate the feature maps. In our case, the features reaching the output layer are the local PDE discretization. The time integrators are then used to evolve the networks output over time for $k$ steps forward and backward, which is then used to compute the loss.}
    \label{workflow}
\end{figure}

\begin{prop}[\textbf{A rational function approximates solution-adaptive discretizations}]\label{prop3}
The solution-adaptive WENO approximation of a function $f$ at points $x_{i \pm 1/2}$ generates stencil weights $\nu$ that are rational functions $\mathcal{G}$ in the stencil values $\mathbf{f}_m (x_i)$, i.e., 
\begin{align}
    \nu_j^{\pm} &= \mathcal{G}\big( \textbf{f}_m(x_i) \big) = \frac{\textbf{p}(\mathbf{f}_m^{\pm}(x_i))}{\textbf{q}(\mathbf{f}_m^{\pm}(x_i))}, \textrm{ where } \mathbf{f}_m^{\pm}(x_i) = \{f(x_j): x_j \in S_m^{\pm}(x_i) \}.
\end{align}
Here, $\textbf{p} : \mathbb{R}^{2m+1} \rightarrow \mathbb{R}$,  $\textbf{q} : \mathbb{R}^{2m+1} \rightarrow \mathbb{R}$ are polynomials of degree  $\geq 1$. Numerical corrections are usually applied to prevent the denominator polynomial $\textbf{q}$ from going to zero (\cite{osher2004level}). The degree of the rational function $\mathcal{G} \left( \cdot \right)$ is the maximum of the degrees of its numerator and denominator polynomials, and the denominator polynomial is always strictly positive. (Proof in Appendix~\ref{WENO_appendix}).
\end{prop}

\begin{coro}[\textbf{A rational function approximates solution-adaptive derivative discretizations}] \label{coro1}
From Proposition (\ref{prop3}) for the WENO-approximated function $\hat{f}$, it follows that its symmetric difference 
\begin{equation}
   \left. \frac{\partial f}{\partial x} \right\vert_{x=x_i} = \frac{\hat{f}_{i+1/2} - \hat{f}_{i-1/2}}{\Delta x} + \mathcal{O}\left( \Delta x^r\right)
\end{equation}
is also a rational function in the stencil values $\mathbf{f}_m (x_i)$
(\cite{shu1998essentially}).
\end{coro}
In finite-difference methods, the function $f$ is a polynomial flux, i.e., $f(u) = c_1 u^2 + c_2 u$, whereas in finite-volume methods, the function $f$ is the cell average $\Bar{u} = \frac{1}{\Delta x} \!\int_{x-\Delta x/2}^{x+\Delta x/2}   u(\xi) d\xi $.

\begin{coro} [\textbf{Rational functions approximate solution-adaptive PDE discretizations}] \label{coro2}
From Proposition (\ref{prop2}) and Corollary (\ref{coro1}), a WENO-like local solution-adaptive discretization of any term in the continuous PDE right-hand side $\mathcal{N} \left( \cdot \right)$ can be approximated by a rational function in the stencil values $\mathbf{u}_m(x_i)$, up to any desired order of accuracy $r$ on a grid with resolution $\Delta x$, even if the solution $u$ admits discontinuities:
\begin{equation}
    \mathcal{N} = \mathcal{N}_d \big( \mathbf{u}_m(x_i), \Xi \big) + \mathcal{O}(\Delta x^r), \textrm{ where } \mathbf{u}_m(x_i) = \{ u(x_j): x_j \in S_m(x_i) \},
\end{equation}
where $\mathcal{N}_d (\cdot)$ is a rational function and also the discrete approximation to the right-hand side of the PDE. The discretization is usually problem specific, hence the dependence on PDE coefficients $\Xi$.
\end{coro}

\section{Connecting MLP convolutional (mlpconv) layers, rational functions, and solution-adaptive PDE operators} \label{section4}
In CNN architectures, the $k^{th}$ feature map channel for a 1D input patch is computed as $f_{i}^{k} = \sigma \left( \omega^{k} x_i\right)$, where $i$ is the pixel index,  $x_i$ is the input patch centered at location $i$, and $\omega$ is the filter moved across the input patch. The computation performed by a single mlpconv layer is:
\begin{equation}\label{mlp_compute}
    \mathbf{y}_q^1 = \sigma \left( \mathbf{W}^1 \mathbf{x} + \mathbf{b}_1\right), \hdots \hdots \mathbf{y}_q^{n_l} = \sigma \left( \mathbf{W}^{n_l} \mathbf{y}_q^{n_l-1} + \mathbf{b}_n\right),
\end{equation}
\begin{wrapfigure}{r}{0.4\textwidth} 
\vspace{-1.0em}
\centering
\includegraphics[width=0.4\textwidth]{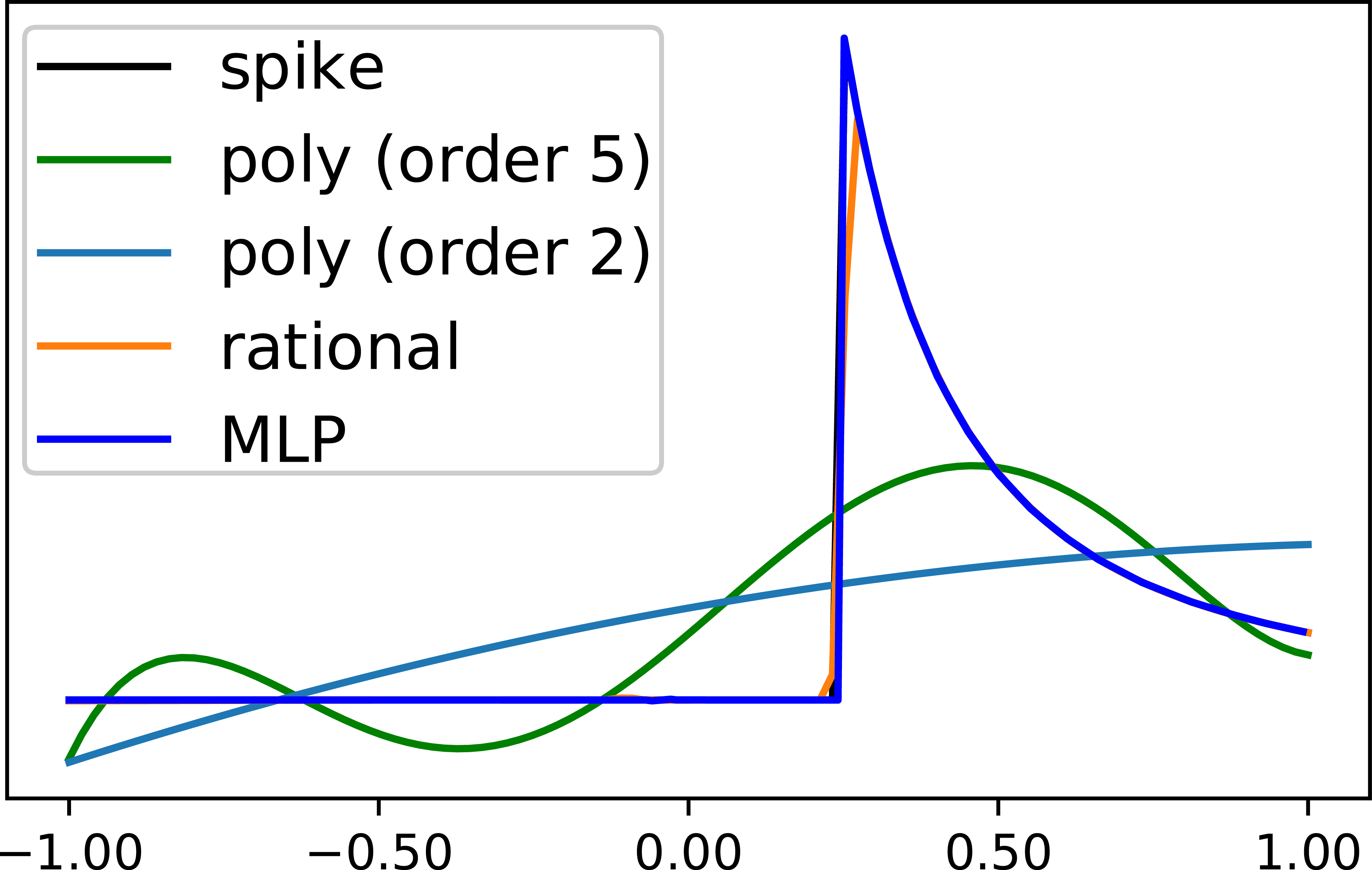}
\caption{\footnotesize Rational, polynomial, and MLP fit to the spike function $f(x) = \left(a + a\vert x - b\vert/(x-b) \right) (1/(x+c)^2)$ with $a = 0.5, b = 0.25$, and $c=0.1$. The rational function fit uses Newman polynomials to approximate $\vert x\vert$.}
\label{rational}
\vspace{-2em}
\end{wrapfigure}
where $\theta = \{ \textbf{W}_q, \textbf{b}_q \}_{q=1,2,\ldots , n_l}$ are the trainable weights and biases, respectively, and $\sigma$ is the (nonlinear) activation function.
The mlpconv layer can thus be seen as cascaded cross-channel parametric pooling on a normal convolutional layer (\cite{lin2013network}), which allows complex and learnable interactions across channels for better abstraction of the input data. An important property of MLPs, which we exploit in this work, is their effectiveness in approximating rational functions, and vice-versa (\cite{telgarsky2017neural}), as stated in the following proposition.

\begin{prop} [\textbf{MLPs $\boldsymbol{\epsilon}$-approximate rational functions --- \cite{telgarsky2017neural}}] \label{prop4} For polynomials $\textbf{p}: [0,1]^d \rightarrow [-1,+1]$ and $\textbf{q}: [0,1]^d \rightarrow [2^{-k},1]$  of degree $\leq r$, each with $\leq s$ monomials, there exists a fully connected Rectified Linear Unit (ReLU) network $f:[0,1]^d \rightarrow \mathbb{R}$ such that
\begin{equation}
    \sup_{ x \in [0,1]^d} \Bigg \vert f(x) - \frac{\textbf{p}(x)}{\textbf{g}(x)} \Bigg \vert \leq \epsilon ,
\end{equation}
for $\epsilon \in (0,1]$. The reverse of the above statement is also true (see proof in \cite{telgarsky2017neural} (3.1)).
\end{prop}

This enables MLPs to generate features that mimic WENO-like solution-adaptive discretizations of PDEs. The rational-function nature of WENO schemes originates from the need to compute reciprocals of smoothness indicators (see Appendix \ref{WENO_appendix}). This requires that any surrogate  model of solution-adaptive discretizations also  approximates division (reciprocal) operations. It is clear from Eqs.~(\ref{linear_derv}) and (\ref{non_linear_convolutions}) that the class of polynomials (i.e., nonlinear convolutions) cannot approximate division operation. Thus, methods based on multiplication of filters (\cite{long2019pde}) are not solution-adaptive. Interestingly, \cite{telgarsky2017neural} also showed that both rational functions and MLPs produce higher-fidelity approximations to reciprocal operations, absolute value functions, and ``switch statements'' than polynomials. We illustrate this behavior by comparing the polynomial, rational, and MLP fits to the spike function in Fig.~\ref{rational}. The rational and MLP approximations both closely follow the curve, whereas the polynomial class fails to capture the ``switch statement'' or the division operations that are quintessential for resolving sharp functions. This also explains the results, shown in Fig.~\ref{steep}, where WENO and STENCIL-NET (MLP-based) methods readily resolve discontinuities in the solution of advection equation. Remarkably, STENCIL-NET can advect sharp features on $4\times$ coarser grids than the WENO scheme (Fig.\ref{steep} D).

Based on Corollary (\ref{coro2}) and Proposition (\ref{prop4}), we posit that WENO-like solution-adaptive discretization of PDE operators can be adequately achieved via parametric pooling with mlpconv layers. By sliding an MLP over the reconstructed input state-variable vector $\mathbf{U}^{n}$ (see Fig.~\ref{workflow}), we can map the input on the local stencil patch to complex discretization features. 

\begin{coro} [\textbf{STENCIL-NET MLP $\boldsymbol{\epsilon}$-approximates solution-adaptive PDE discretization}] \label{stencilnet_coro} A direct consequence of Corollary (\ref{coro2}) and Proposition (\ref{prop4}) is that a fully connected ReLU network with stencil values $\mathbf{u}_m(x_i)$ as input can approximate the discrete nonlinear function $\mathcal{N}_d(\cdot)$ at point $x_i$, even if the solution $u$ admits discontinuities across the stencil:
\begin{equation}
\Big \vert \hat{\mathcal{N}}_{\theta} \big( \mathbf{u}_m(x_i), \Xi \big)  - \mathcal{N}_d \big( \mathbf{u}_m(x_i) , \Xi \big) \Big \vert \leq \epsilon.
\end{equation}
Here, $\hat{\mathcal{N}}_{\theta}: \mathbb{R}^{2m+1} \rightarrow \mathbb{R}$ and $\epsilon \in (0,1]$.
\end{coro}
The grid spacing $\Delta x$ is an implicit input to the network, given that sampling happens on a grid with spacing $\Delta x$. Thus, the learned STENCIL-NET parameterizations are specific to a given grid resolution $\Delta x$.

\section{Numerical experiments}\label{section5}
\begin{wrapfigure}{r}{0.5\textwidth}     
\vspace{-3.0em}
\centering
\includegraphics[width=0.48\textwidth]{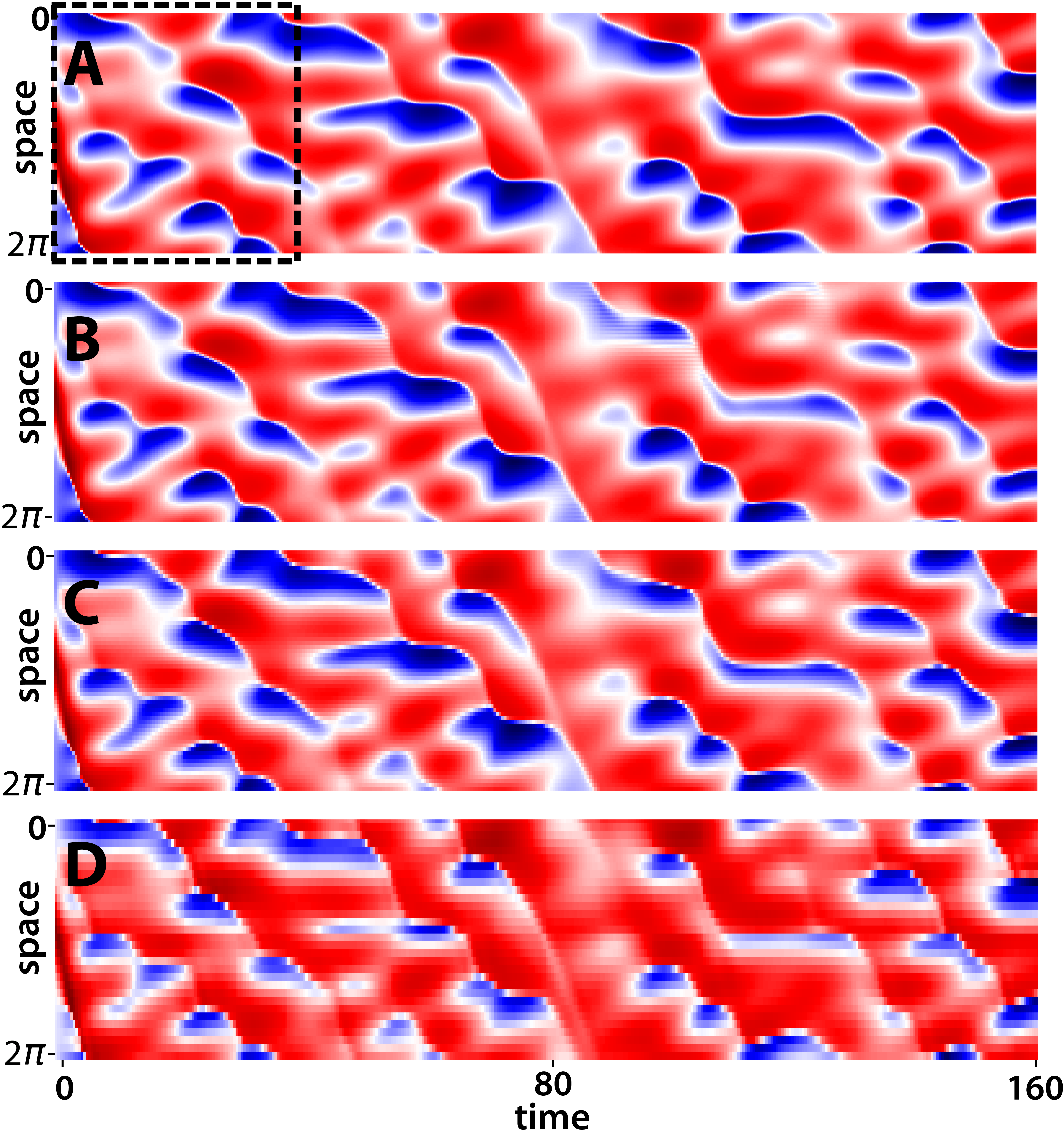}
    \caption{\small \textbf{Forced Burgers' prediction with STENCIL-NET}: (A) Fifth-order WENO solution with $N_x = 256$. (B,C,D) Predictions of STENCIL-NET on ($2\times, 4\times, 8\times $) coarser grids, respectively. The dashed box in A is the domain used for training.}
    \label{burgers_prediction}
    \vspace{-3.5em}
\end{wrapfigure}

We apply STENCIL-NET to learn solution-adaptive discretizations of nonlinear PDEs. For all problems discussed here, we use a single mlpconv unit with 3 hidden layers, each with 64 nodes and Exponential Linear Unit (ELU) nonlinearity. The input to the network is the solution $\mathbf{u}_m$ on the stencil with $m=3$ in all examples. The schematic of the network architecture is shown in Fig.~\ref{workflow}. Details regarding model training are discussed in Appendix \ref{PDE_optim_appendix},
and detailed descriptions of the mathematical models and numerical methods used to generate the data presented in the Appendix \ref{data_gen_append}.

In the following sections, we present numerical experiments that demonstrate three distinct applications of STENCIL-NET based on learning local solution-adaptive discretization.

\subsection{STENCIL-NET: A fast and accurate predictive surrogate} \label{burgers_section}
The forced Burgers' equation with a nonlinear forcing term can produce rich and sharply varying solutions. The forcing term also introduces randomness that can help explore the solution manifold (\cite{bar2019learning}). This presents an ideal challenge for testing STENCIL-NET's discretization capabilities. We train STENCIL-NET on fifth-order WENO solutions of the forced Burgers' equation at different spatial resolutions $\Delta x_c = (C\Delta x)$, where $C$ is the sub-sampling factor, and $\Delta x$ the resolution of the fine-grid solution (Fig.~\ref{burgers_prediction}). More details on the test case are in Appendix~\ref{data_gen_append}.
\begin{wrapfigure}{l}{0.5\textwidth}
\vspace{0em}
\centering
\includegraphics[width=0.45\textwidth]{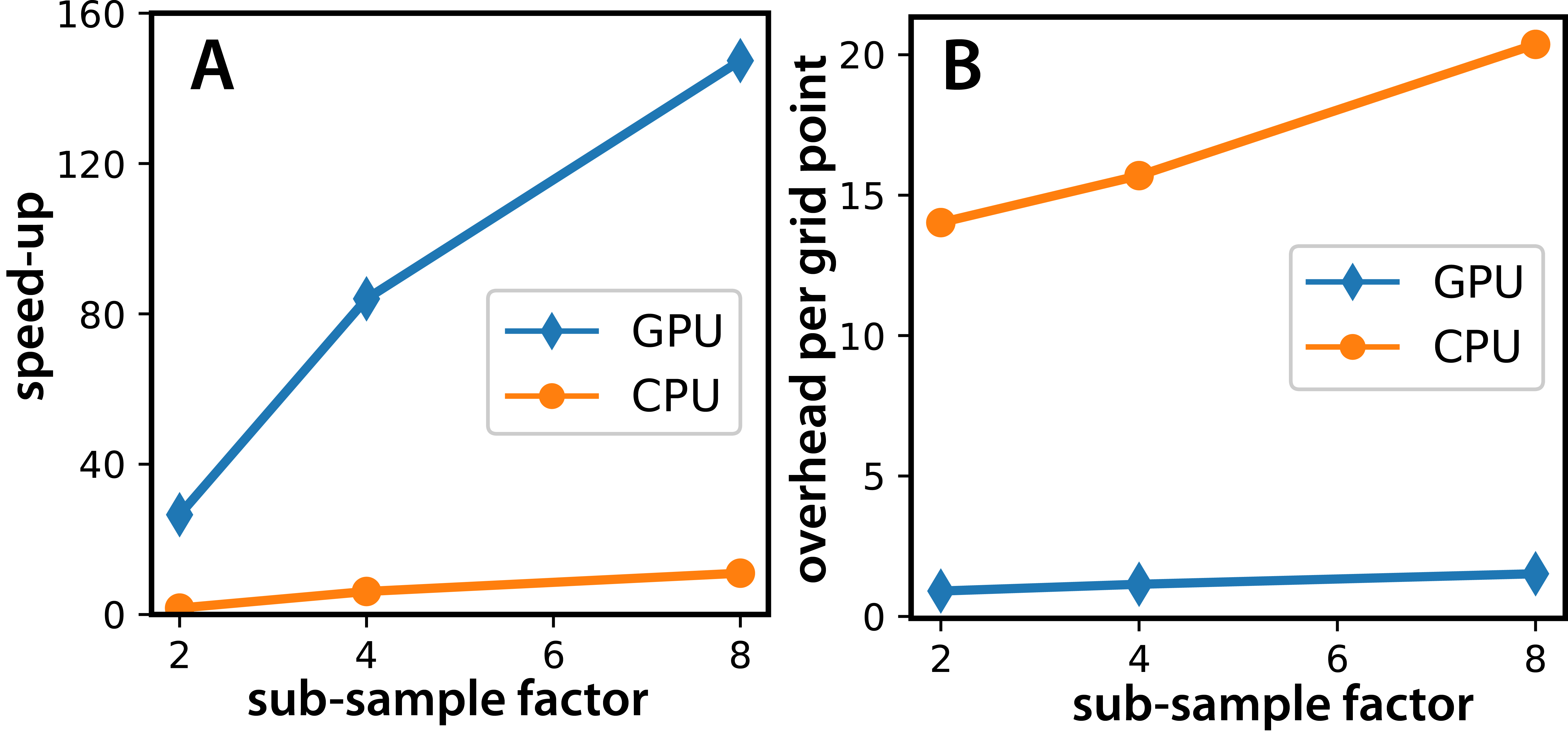}
    \caption{\small \textbf{STENCIL-NET for data-driven acceleration}: Computational speed-ups (A) and overhead per grid point (B, see Appendix \ref{speedup_appendix}) on CPUs and GPUs for differently sub-sampled forced Burgers' cases.}
    \label{speedup}
\end{wrapfigure}

\begin{wrapfigure}{r}{0.45\textwidth}
\vspace{-16em}
\centering
\includegraphics[width=0.42\textwidth]{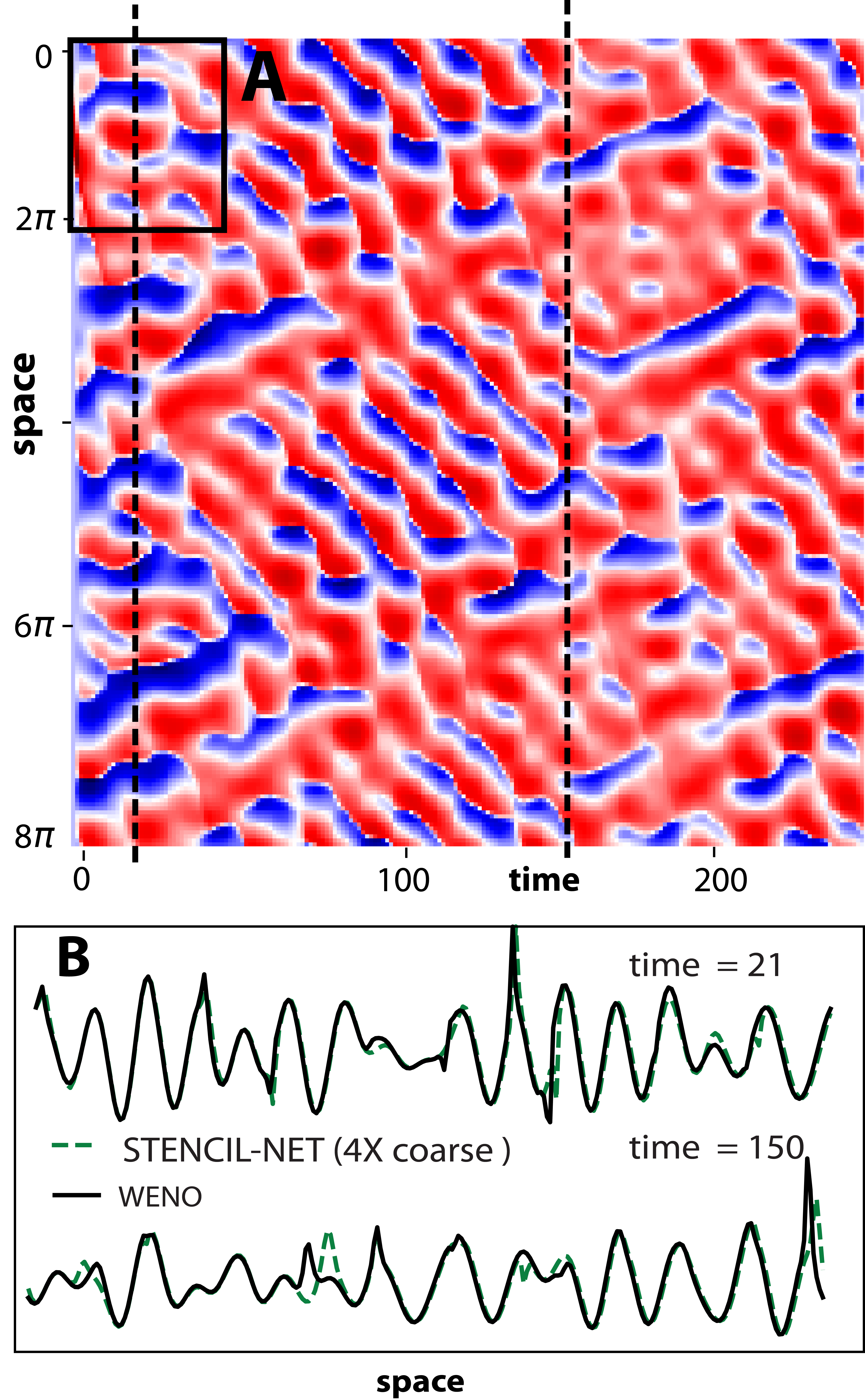}
    \caption{\small \textbf{STENCIL-NET solution on larger spatial domains for longer times} (A): STENCIL-NET predictions on $4\times$ coarser grid. (B): Comparison between WENO evalutation of the discrete operator $\mathcal{N}_d$ and mlpconv layer output $\hat{\mathcal{N}}_{\theta}$ at $21$s and $150$s (dashed lines in A). The STENCIL-NET was only trained on the solution in the domain marked by the solid rectangle in A. }
    \label{burgers_prediction_large}
    \vspace{-1.5em}
\end{wrapfigure}

In Fig.~\ref{burgers_prediction}, we show a comparison between fifth-order accurate WENO solution on a fine grid $(\Delta x = L/N_x, N_x = 256)$ with domain length $L=2\pi$ and the corresponding STENCIL-NET predictions on coarsened grids with $\Delta x_c = C\Delta x$, for sub-sampling factors $C = \{2,4,8\}$. STENCIL-NET produces \textit{stable} solutions on coarser spatio-temporal grids for longer times with little compromise on accuracy.
In contrast to \cite{kim2019deep}, which uses all volumetric time frames as input to the network, STENCIL-NET exploits parametric pooling on local stencil patches, which enables extrapolation to larger spatial domains and longer times with negligible loss in accuracy (Fig.~\ref{burgers_prediction_large}A,B). More details on accuracy and convergence are presented in Appendix \ref{more_results_appendix}.

In Fig.~\ref{speedup}, we show performance metrics quantifying the speed-up achieved by STENCIL-NET predictions over the base-line WENO solvers on large domains $L = 64\pi$ discretized with $N_x = 8192$ grid points. The speed-up for STENCIL-NET is mainly due to two factors: 1) Reducing the number of grid points by a factor of at least $C^{d+2}$, where $d$ is the spatial dimension; 2) Substituting low compute-to-memory operations (classical numerics) with dense matrix multiplication reaching near-peak performance on both CPU and GPU architectures (\cite{jouppi2017datacenter,zhuang2020learned}). We find speed-ups between $25\ldots 150\times$ on GPUs and $2\ldots 14 \times$ on CPUs for different sub-sampling factors shown in Fig.~\ref{speedup}A. The base-line numerical solver was a fully vectorized WENO solver for Burgers' equation. We also quantify the overhead per grid point (Fig.~\ref{speedup}B) from substituting WENO operators with the forward-pass evaluation (Eq.~\ref{mlp_compute}) of the mlpconv layer. All GPU runtime estimates were measured on Nvidia GeForce GTX 1080 accelerator. Under the fair assumption of a linear relation between computational overhead and dimension $d$, we argue that the potential speed-up from STENCIL-NET coarse-graining increases exponentially with the problem's spatial dimension. We provide an estimate of the achievable speed-up in Appendix \ref{speedup_appendix}.

\subsection{STENCIL-NET: An autonomous predictor of chaotic dynamics} \label{KS_section}
\begin{figure}[!htbp]
    \centering
    \vspace{-0.5em}
    \includegraphics[width=0.95\textwidth]{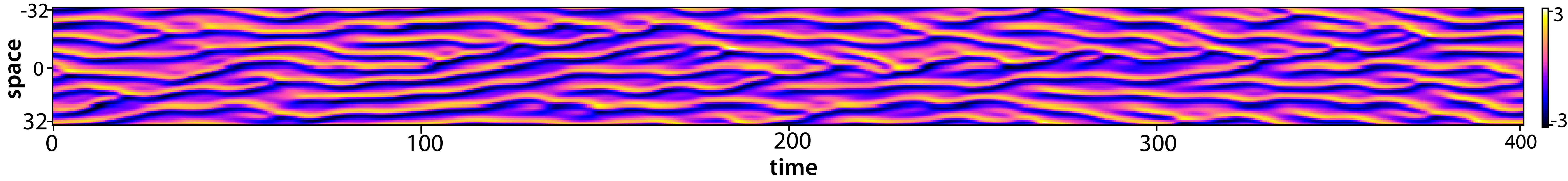}
    \vspace{-0.5em}
    \caption{\small \textbf{Equation-free prediction of chaotic spatio-temporal dynamics:} STENCIL-NET run as autonomous system for predicting long-term chaotic dynamics on a $4\times$ coarsened grid.}
    \label{autonomous}
\end{figure}

Nonlinear dynamical systems, like the Kuramoto-Sivashinsky (KS) model, can describe chaotic spatio-temporal dynamics. The KS equation exhibits varying levels of chaos depending on the bifurcation parameter $L$ (domain size) of the system (\cite{vlachas2018data}). Given their sensitivity to numerical errors, KS equations are usually solved using spectral methods. 

\begin{wrapfigure}{r}{0.5\textwidth}
\centering
\vspace{-0.5em}
\includegraphics[width=0.48\textwidth]{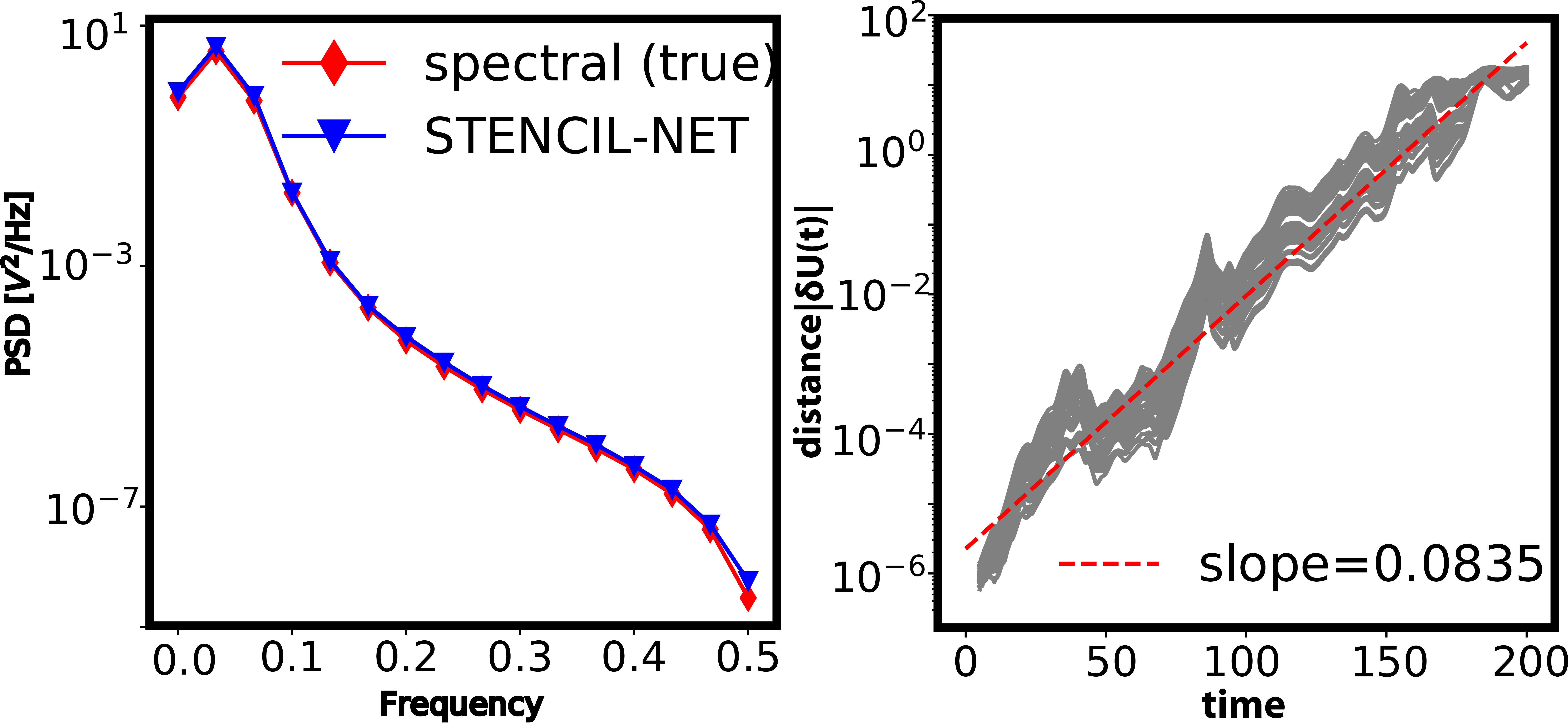}
    \caption{\small \textbf{Quantitative metrics for chaotic systems:} (Left) Power spectrum comparison between STENCIL-NET and ground-truth spectral solution. (Right) Growth of distance between nearby trajectories characterizing the maximum Lyaponuv exponent (slope of dashed line), compared to ground truth $\approx 0.084$ (\cite{edson2019lyapunov}).}
    \label{KS_metrics}
    \vspace{-1em}
\end{wrapfigure}

However, recent data-driven models using Recurrent Neural Networks (RNNs) (\cite{patankar2018numerical,vlachas2018data,vlachas2020backpropagation}) have shown success in predicting chaotic systems. This is mainly because RNNs are able to capture long-term temporal correlations and implicitly identify the required embedding for forecasting.

We challenge these results using STENCIL-NET for long-term stable prediction of chaotic dynamics. The training data are from spectral solutions of the KS equation for domain size $L=64$ (see Appendix \ref{data_gen_append}). In Fig.~\ref{autonomous}, we show STENCIL-NET predictions on $4\times$ coarser grid with different initial condition and for longer time runs ($ 8\times$) than available during training. We find that STENCIL-NET on $4\times$ coarser grid accurately captures the long-term spectral statistics of the nonlinear dynamical system (Fig.~\ref{KS_metrics}, left) and produces an accurate estimate of the maximum Lyapunov exponent (Fig.~\ref{KS_metrics}, right).

\subsection{STENCIL-NET: A de-noising method} \label{KDv_section}
\begin{figure}[!htbp]
    \centering
    \includegraphics[width=1.0\textwidth]{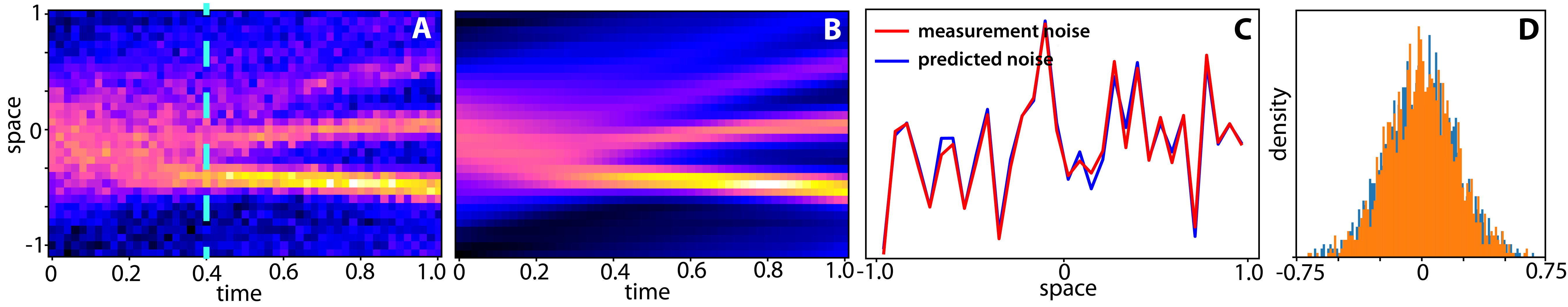}
    \vspace{-1.2em}
    \caption{\small \textbf{STENCIL-NET for de-noising dynamical data:} (A): KdV solution with additive noise used for training. (B): STENCIL-NET de-noised prediction. (C): Comparison between STENCIL-NET noise estimates and actual noise at $t=0.4$\,s. (D) Comparison between STENCIL-NET's estimated noise distribution (orange) and true noise distribution (blue).}
    \label{denoise}
    \vspace{-1.2em}
\end{figure}

The discrete time-stepping constraints force any STENCIL-NET prediction to follow a smooth time trajectory. This property can be exploited for filtering true dynamics from noise. We demonstrate this de-noising capability using numerical solution data for the Korteweg-de Vries (KdV) equation with artificial additive noise (details in Appendix \ref{data_gen_append}). From the noisy data, STENCIL-NET is able to learn an accurate and consistent discretization of the KdV equation, enabling it to decompose dynamics and noise. Moreover, the network is able to quantitatively reconstruct point-wise estimates of the noise (Fig.~\ref{denoise}C) as well as the noise distribution (Fig.~\ref{denoise}D) without any physics prior or model of the noise distribution. The detailed description of the optimization formulation for de-noising can be found in Appendix \ref{PDE_optim_appendix}.

\section{Conclusion}
We have presented the STENCIL-NET architecture for data-driven learning of solution-adaptive discretizations of nonlinear PDEs. STENCIL-NET uses patch-wise parametric pooling and discrete time integration constraints to learn coarse-grained numerical predictions directly from spatio-temporal data, without knowledge of the underlying symbolic PDE. We have shown that STENCIL-NET can be used as fast and accurate predictive surrogate for accelerating numerical simulations, for model-free forecasting of chaotic dynamics, and for measurement de-noising. The design of the STENCIL-NET architecture rests on a formal connection between MLPs, rational functions, and solution-adaptive discretization schemes for PDEs, which we have presented here. 

Future generalizations of STENCIL-NET include an extension to 2D and 3D problems, multi-resolution numerical schemes, and mesh-free particle methods. We hope that the present results provide motivation for future works and for including STENCIL-NET into existing simulations frameworks.

\section*{Broader Impact}

Often in science and engineering, measurement data from a dynamical processes in space and time are available, but a first-principles mathematical model is either difficult to formulate or too simplistic to realistically capture the data. Numerical simulation methods can then not be used to predict system behavior. This problem is particularly prevalent in areas such as biology, medicine, environmental science, economy, and finance. 

STENCIL-NET fills this niche by providing a purely data-driven, equation-free way of learning predictive computable models of complex nonlinear space-time dynamics. It can be used to produce rapid coarse-grained predictions beyond the space and time horizon it was trained for, even for chaotic dynamics, and for accurate de-noising of measurement data. Beyond the applications shown here, we anticipate the fast and accurate STENCIL-NET predictions to be useful in areas including real-time computational steering, computer graphics, phase-space exploration, uncertainty quantification, and computer vision.

\begin{ack}
This work was supported by the German Research Foundation (DFG) as part of the Federal Cluster of Excellence ``Physics of Life'' under code EXC 2068, and by the German Federal Ministry of Education and Research (BMBF) as part of the National Center for Scalable Data Analytics and Artificial Intelligence (ScaDS.AI).
\end{ack}

\medskip

\small

\bibliography{ref}

\begin{thebibliography}{36}
\providecommand{\natexlab}[1]{#1}
\providecommand{\url}[1]{\texttt{#1}}
\expandafter\ifx\csname urlstyle\endcsname\relax
  \providecommand{\doi}[1]{doi: #1}\else
  \providecommand{\doi}{doi: \begingroup \urlstyle{rm}\Url}\fi

\bibitem[Patankar(2018)]{patankar2018numerical}
Suhas Patankar.
\newblock \emph{Numerical heat transfer and fluid flow}.
\newblock Taylor \& Francis, 2018.

\bibitem[Moin(2010)]{moin2010fundamentals}
Parviz Moin.
\newblock \emph{Fundamentals of engineering numerical analysis}.
\newblock Cambridge University Press, 2010.

\bibitem[Osher et~al.(2004)Osher, Fedkiw, and Piechor]{osher2004level}
Stanley Osher, Ronald Fedkiw, and K~Piechor.
\newblock Level set methods and dynamic implicit surfaces.
\newblock \emph{Appl. Mech. Rev.}, 57\penalty0 (3):\penalty0 B15--B15, 2004.

\bibitem[Shu(1998)]{shu1998essentially}
Chi-Wang Shu.
\newblock Essentially non-oscillatory and weighted essentially non-oscillatory
  schemes for hyperbolic conservation laws.
\newblock In \emph{Advanced numerical approximation of nonlinear hyperbolic
  equations}, pages 325--432. Springer, 1998.

\bibitem[Bijl et~al.(2013)Bijl, Lucor, Mishra, and Schwab]{bijl2013uncertainty}
Hester Bijl, Didier Lucor, Siddhartha Mishra, and Christoph Schwab.
\newblock \emph{Uncertainty quantification in computational fluid dynamics},
  volume~92.
\newblock Springer Science \& Business Media, 2013.

\bibitem[Mishra(2018)]{mishra2018machine}
Siddhartha Mishra.
\newblock A machine learning framework for data driven acceleration of
  computations of differential equations.
\newblock \emph{arXiv preprint arXiv:1807.09519}, 2018.

\bibitem[Tompson et~al.(2017)Tompson, Schlachter, Sprechmann, and
  Perlin]{tompson2017accelerating}
Jonathan Tompson, Kristofer Schlachter, Pablo Sprechmann, and Ken Perlin.
\newblock Accelerating eulerian fluid simulation with convolutional networks.
\newblock In \emph{Proceedings of the 34th International Conference on Machine
  Learning-Volume 70}, pages 3424--3433. JMLR. org, 2017.

\bibitem[Kim et~al.(2019)Kim, Azevedo, Thuerey, Kim, Gross, and
  Solenthaler]{kim2019deep}
Byungsoo Kim, Vinicius~C Azevedo, Nils Thuerey, Theodore Kim, Markus Gross, and
  Barbara Solenthaler.
\newblock Deep fluids: A generative network for parameterized fluid
  simulations.
\newblock In \emph{Computer Graphics Forum}, volume~38, pages 59--70. Wiley
  Online Library, 2019.

\bibitem[Pathak et~al.(2017)Pathak, Lu, Hunt, Girvan, and Ott]{pathak2017using}
Jaideep Pathak, Zhixin Lu, Brian~R Hunt, Michelle Girvan, and Edward Ott.
\newblock Using machine learning to replicate chaotic attractors and calculate
  lyapunov exponents from data.
\newblock \emph{Chaos: An Interdisciplinary Journal of Nonlinear Science},
  27\penalty0 (12):\penalty0 121102, 2017.

\bibitem[Vlachas et~al.(2018)Vlachas, Byeon, Wan, Sapsis, and
  Koumoutsakos]{vlachas2018data}
Pantelis~R Vlachas, Wonmin Byeon, Zhong~Y Wan, Themistoklis~P Sapsis, and
  Petros Koumoutsakos.
\newblock Data-driven forecasting of high-dimensional chaotic systems with long
  short-term memory networks.
\newblock \emph{Proceedings of the Royal Society A: Mathematical, Physical and
  Engineering Sciences}, 474\penalty0 (2213):\penalty0 20170844, 2018.

\bibitem[Lin et~al.(2013)Lin, Chen, and Yan]{lin2013network}
Min Lin, Qiang Chen, and Shuicheng Yan.
\newblock Network in network.
\newblock \emph{arXiv preprint arXiv:1312.4400}, 2013.

\bibitem[Brunton et~al.(2016)Brunton, Proctor, and
  Kutz]{brunton2016discovering}
Steven~L Brunton, Joshua~L Proctor, and J~Nathan Kutz.
\newblock Discovering governing equations from data by sparse identification of
  nonlinear dynamical systems.
\newblock \emph{Proceedings of the national academy of sciences}, 113\penalty0
  (15):\penalty0 3932--3937, 2016.

\bibitem[Rudy et~al.(2017)Rudy, Brunton, Proctor, and Kutz]{rudy2017data}
Samuel~H Rudy, Steven~L Brunton, Joshua~L Proctor, and J~Nathan Kutz.
\newblock Data-driven discovery of partial differential equations.
\newblock \emph{Science Advances}, 3\penalty0 (4):\penalty0 e1602614, 2017.

\bibitem[Maddu et~al.(2019)Maddu, Cheeseman, Sbalzarini, and
  M{\"u}ller]{maddu2019stability}
Suryanarayana Maddu, Bevan~L Cheeseman, Ivo~F Sbalzarini, and Christian~L
  M{\"u}ller.
\newblock Stability selection enables robust learning of partial differential
  equations from limited noisy data.
\newblock \emph{arXiv preprint arXiv:1907.07810}, 2019.

\bibitem[Raissi et~al.(2017)Raissi, Perdikaris, and
  Karniadakis]{raissi2017physics}
Maziar Raissi, Paris Perdikaris, and George~Em Karniadakis.
\newblock Physics informed deep learning (part {I}): Data-driven solutions of
  nonlinear partial differential equations.
\newblock \emph{arXiv preprint arXiv:1711.10561}, 2017.

\bibitem[Raissi et~al.(2020)Raissi, Yazdani, and Karniadakis]{raissi2020hidden}
Maziar Raissi, Alireza Yazdani, and George~Em Karniadakis.
\newblock Hidden fluid mechanics: Learning velocity and pressure fields from
  flow visualizations.
\newblock \emph{Science}, 367\penalty0 (6481):\penalty0 1026--1030, 2020.

\bibitem[Raissi et~al.(2018)Raissi, Perdikaris, and
  Karniadakis]{raissi2018multistep}
Maziar Raissi, Paris Perdikaris, and George~Em Karniadakis.
\newblock Multistep neural networks for data-driven discovery of nonlinear
  dynamical systems.
\newblock \emph{arXiv preprint arXiv:1801.01236}, 2018.

\bibitem[Rudy et~al.(2019)Rudy, Kutz, and Brunton]{rudy2019deep}
Samuel~H Rudy, J~Nathan Kutz, and Steven~L Brunton.
\newblock Deep learning of dynamics and signal-noise decomposition with
  time-stepping constraints.
\newblock \emph{Journal of Computational Physics}, 396:\penalty0 483--506,
  2019.

\bibitem[Long et~al.(2017)Long, Lu, Ma, and Dong]{long2017pde}
Zichao Long, Yiping Lu, Xianzhong Ma, and Bin Dong.
\newblock {PDE}-{N}et: Learning {PDE}s from data.
\newblock \emph{arXiv preprint arXiv:1710.09668}, 2017.

\bibitem[Long et~al.(2019)Long, Lu, and Dong]{long2019pde}
Zichao Long, Yiping Lu, and Bin Dong.
\newblock {PDE-N}et 2.0: Learning {PDE}s from data with a numeric-symbolic
  hybrid deep network.
\newblock \emph{Journal of Computational Physics}, 399:\penalty0 108925, 2019.

\bibitem[Bar-Sinai et~al.(2019)Bar-Sinai, Hoyer, Hickey, and
  Brenner]{bar2019learning}
Yohai Bar-Sinai, Stephan Hoyer, Jason Hickey, and Michael~P Brenner.
\newblock Learning data-driven discretizations for partial differential
  equations.
\newblock \emph{Proceedings of the National Academy of Sciences}, 116\penalty0
  (31):\penalty0 15344--15349, 2019.

\bibitem[Higdon(2006)]{higdon2006numerical}
Robert~L Higdon.
\newblock Numerical modelling of ocean circulation.
\newblock \emph{Acta Numerica}, 15:\penalty0 385--470, 2006.

\bibitem[Veynante and Vervisch(2002)]{veynante2002turbulent}
Denis Veynante and Luc Vervisch.
\newblock Turbulent combustion modeling.
\newblock \emph{Progress in energy and combustion science}, 28\penalty0
  (3):\penalty0 193--266, 2002.

\bibitem[Prost et~al.(2015)Prost, J{\"u}licher, and Joanny]{prost2015active}
Jacques Prost, Frank J{\"u}licher, and Jean-Fran{\c{c}}ois Joanny.
\newblock Active gel physics.
\newblock \emph{Nature Physics}, 11\penalty0 (2):\penalty0 111--117, 2015.

\bibitem[Fornberg(1988)]{fornberg1988generation}
Bengt Fornberg.
\newblock Generation of finite difference formulas on arbitrarily spaced grids.
\newblock \emph{Mathematics of computation}, 51\penalty0 (184):\penalty0
  699--706, 1988.

\bibitem[Schrader et~al.(2010)Schrader, Reboux, and
  Sbalzarini]{schrader2010discretization}
Birte Schrader, Sylvain Reboux, and Ivo~F Sbalzarini.
\newblock Discretization correction of general integral {PSE} operators for
  particle methods.
\newblock \emph{Journal of Computational Physics}, 229\penalty0 (11):\penalty0
  4159--4182, 2010.

\bibitem[Zoumpourlis et~al.(2017)Zoumpourlis, Doumanoglou, Vretos, and
  Daras]{zoumpourlis2017non}
Georgios Zoumpourlis, Alexandros Doumanoglou, Nicholas Vretos, and Petros
  Daras.
\newblock Non-linear convolution filters for {CNN}-based learning.
\newblock In \emph{Proceedings of the IEEE International Conference on Computer
  Vision}, pages 4761--4769, 2017.

\bibitem[Courant et~al.(1952)Courant, Isaacson, and Rees]{courant1952solution}
Richard Courant, Eugene Isaacson, and Mina Rees.
\newblock On the solution of nonlinear hyperbolic differential equations by
  finite differences.
\newblock \emph{Communications on pure and applied mathematics}, 5\penalty0
  (3):\penalty0 243--255, 1952.

\bibitem[Sod(1985)]{sod1985numerical}
Gary~A Sod.
\newblock \emph{Numerical methods in fluid dynamics: initial and initial
  boundary-value problems}.
\newblock Cambridge University Press, 1985.

\bibitem[Telgarsky(2017)]{telgarsky2017neural}
Matus Telgarsky.
\newblock Neural networks and rational functions.
\newblock In \emph{Proceedings of the 34th International Conference on Machine
  Learning-Volume 70}, pages 3387--3393. JMLR. org, 2017.

\bibitem[Jouppi et~al.(2017)Jouppi, Young, Patil, Patterson, Agrawal, Bajwa,
  Bates, Bhatia, Boden, Borchers, et~al.]{jouppi2017datacenter}
Norman~P Jouppi, Cliff Young, Nishant Patil, David Patterson, Gaurav Agrawal,
  Raminder Bajwa, Sarah Bates, Suresh Bhatia, Nan Boden, Al~Borchers, et~al.
\newblock In-datacenter performance analysis of a tensor processing unit.
\newblock In \emph{Proceedings of the 44th Annual International Symposium on
  Computer Architecture}, pages 1--12, 2017.

\bibitem[Zhuang et~al.(2020)Zhuang, Kochkov, Bar-Sinai, Brenner, and
  Hoyer]{zhuang2020learned}
Jiawei Zhuang, Dmitrii Kochkov, Yohai Bar-Sinai, Michael~P Brenner, and Stephan
  Hoyer.
\newblock Learned discretizations for passive scalar advection in a 2-d
  turbulent flow.
\newblock \emph{arXiv preprint arXiv:2004.05477}, 2020.

\bibitem[Edson et~al.(2019)Edson, Bunder, Mattner, and
  Roberts]{edson2019lyapunov}
Russell~A Edson, Judith~E Bunder, Trent~W Mattner, and Anthony~J Roberts.
\newblock Lyapunov exponents of the {K}uramoto--{S}ivashinsky {PDE}.
\newblock \emph{The ANZIAM Journal}, 61\penalty0 (3):\penalty0 270--285, 2019.

\bibitem[Vlachas et~al.(2020)Vlachas, Pathak, Hunt, Sapsis, Girvan, Ott, and
  Koumoutsakos]{vlachas2020backpropagation}
PR~Vlachas, J~Pathak, BR~Hunt, TP~Sapsis, M~Girvan, E~Ott, and P~Koumoutsakos.
\newblock Backpropagation algorithms and reservoir computing in recurrent
  neural networks for the forecasting of complex spatiotemporal dynamics.
\newblock \emph{Neural Networks}, 2020.

\bibitem[Driscoll et~al.(2014)Driscoll, Hale, and
  Trefethen]{driscoll2014chebfun}
Tobin~A Driscoll, Nicholas Hale, and Lloyd~N Trefethen.
\newblock Chebfun guide, 2014.

\bibitem[Kingma and Ba(2014)]{kingma2014adam}
Diederik~P Kingma and Jimmy Ba.
\newblock Adam: A method for stochastic optimization.
\newblock \emph{arXiv preprint arXiv:1412.6980}, 2014.

\end{thebibliography}

\newpage

\appendix
\section{WENO approximation as a rational function}\label{WENO_appendix} 
\setcounter{figure}{0} \renewcommand{\thefigure}{A.\arabic{figure}}
\begin{wrapfigure}{r}{0.5\textwidth} 
\vspace{-2.0em}
\centering
\includegraphics[width=0.45\textwidth]{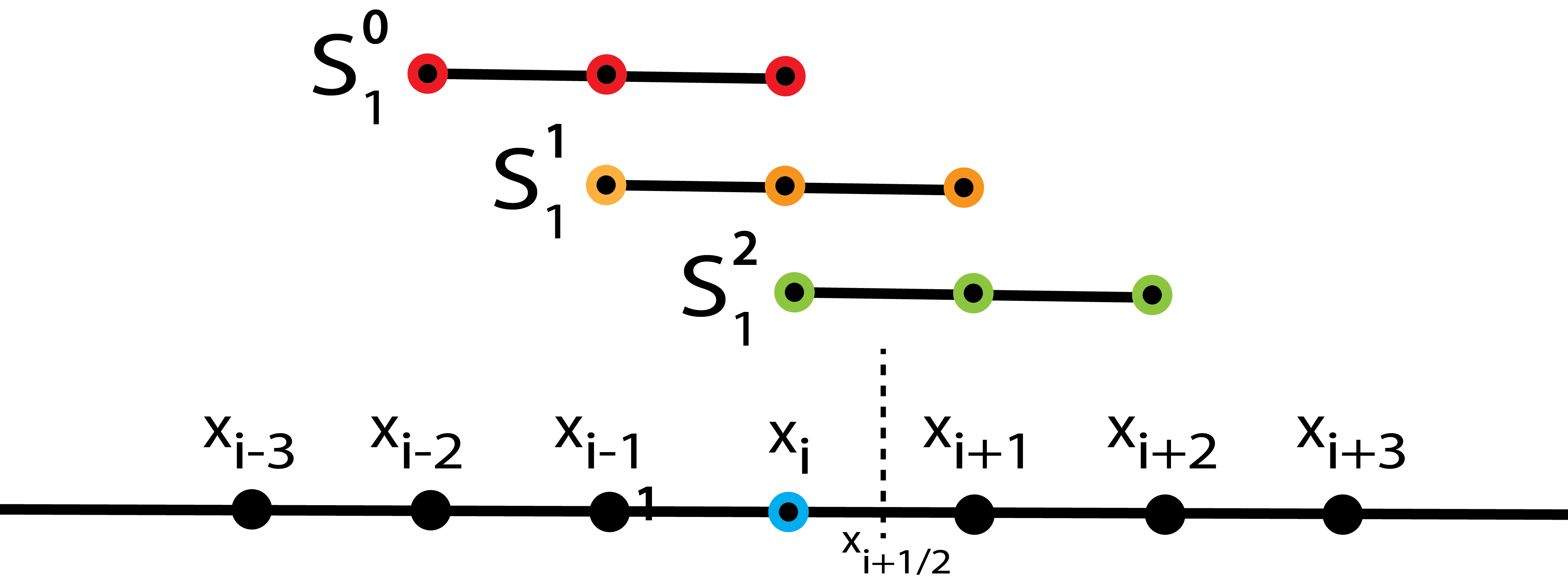}
\caption{\footnotesize Schematic of adaptive stencil choice using fifth-order WENO. The function value at point $x_{i+\frac{1}{2}}$ is computed as the weighted sum of the ENO approximations of each candidate stencil $S_{1}^{k},\, k\in {0,1,2}$.}
\label{stencil_draw_supplementary}
\vspace{-1.0em}
\end{wrapfigure}

ENO and WENO schemes are among the most popular schemes  for piecewise smooth solution approximation across discontinuities. The key idea behind these schemes is to use a nonlinear adaptive procedure to choose the locally smoothest stencil, where discontinuities are avoided to result in smooth and essentially non-oscillatory solutions. 

The WENO approximation of order \mbox{$(2r-1)$} is a weighted sum of the $r$-th order fixed-stencil approximations from $r$ candidate stencils, shown in Fig.~\ref{stencil_draw_supplementary}. This allows to reconstruct the function values $f$ at intermediate point locations, i.e.
\begin{equation} \label{WENO_intermediate}
    f_{i+\frac{1}{2}} = \hat{f}_{i+\frac{1}{2}} + \mathcal{O} (\Delta x^{2r-1}) = \sum_{k=0}^{r-1} \bm{\omega}_k g\left( f_{i+k-r+1},...,f_{i+k} \right) + \mathcal{O} (\Delta x^{2r-1}),
\end{equation}
with the weights $\bm{\omega}_k$ computed in a solution-adaptive fashion. The index $k$ runs through all  candidate stencils used for WENO approximation (see Fig.~\ref{stencil_draw_supplementary}). The function $g$ is the linear map from the function values on the local stencil $S_m^{k}$ to the point of reconstruction:
\begin{equation}\label{linear_map}
    g\left( f_{i+k-r+1},...,f_{i+k} \right) = \hat{f}_{i+\frac{1}{2}}^{k} = \sum_{l=0}^{r-1} c_{k,l}^{r} f_{i+k-r+1+l}.
\end{equation}
The candidate stencil width is chosen such that $\vert S_m^k \vert = 2m+1 = r$. The constant coefficients $c_{k,l}^{r}$ can be computed from the Lagrange interpolation polynomial $\textit{\textbf{p}}_k(x)$ over the $k$-th candidate stencil. This leads to an $r$-th order approximation of the function $f$ using a single candidate stencil $S_m^k$ at the point $x_{i+\frac{1}{2}}$, thus:
\begin{equation} \label{poly_map}
    f_{i+\frac{1}{2}} = \hat{f}_{i+\frac{1}{2}}^{k} + \mathcal{O}(\Delta x^r) = \textit{\textbf{p}}_k(x_{i+\frac{1}{2}}) + \mathcal{O}(\Delta x^r).
\end{equation}

The usual template for the weights $\bm{\omega}_k$ to achieve the essentially non-oscillatory property is:
\begin{equation*}
    \bm{\omega}_k = \frac{\alpha_k}{\alpha_0 + \cdot \cdot \cdot + \alpha_{r-1}}, \quad \alpha_k = \frac{\mathcal{C}_k^r}{ (\varepsilon + \mathbf{IS}_k)^p},
\end{equation*}
where $\mathbf{IS}_k$ is the smoothness indicator of the candidate stencil $S^{k}$, and $p$ is any positive integer. The coefficients $\mathcal{C}_k^r$ are the optimal weights pre-computed in smooth regions. For more details on refer to (\cite{shu1998essentially}). Combining the above relations, we can write the expanded form of the adaptive weights $\bm{\omega}_k$ as,
\begin{equation}\label{complex_weights}
    \bm{\omega}_k =  \frac{\mathcal{C}_k^r \prod_{l=0, l\neq k}^{r-1} (\varepsilon + \mathbf{IS}_l)^p }{\sum_{l=0}^{r-1} C_l^r \prod_{q=0,q\neq l}^{r-1}(\varepsilon + \mathbf{IS}_q)^p }.
\end{equation}

\begin{defn} [\textbf{Smoothness indicators}] \label{def_SI}
For each candidate stencil $S_m^k(x_i)$ we can construct a polynomial $\textbf{p}_k$ of degree $2m$ interpolating the stencil values $\mathbf{f}_m$. The smoothness indicator of the candidate stencil $S_m^k$ is then defined as:
\begin{equation}\label{SI}
    \mathbf{IS}_k =  \sum_{l=1}^{2m} (\Delta x)^{2l-1} \displaystyle\int_{x_{i-1/2}}^{x_{i+1/2}} \left( \frac{\partial^l \textbf{p}_k}{\partial x^l} \right)^{\!\! 2} dx.
\end{equation}
\end{defn}
The term $(\Delta x)^{2l-1}$ removes the $\Delta x$-dependence (\cite{shu1998essentially}). The right-hand side of Eq.~\ref{SI} is the sum of the $L^2$ norms of all derivatives of the polynomial $\textit{\textbf{p}}_k$ interpolating over the stencil $S_m^k$. Using Eqs.~\ref{linear_map} and \ref{poly_map}, we can formulate the following proposition:
\begin{prop} [\textbf{Discrete smoothness indicators are polynomials in the stencil values}] \label{discrete_SI} 
From definition \ref{def_SI} and proposition \ref{prop2} we readily observe that
\begin{equation}
    \mathbf{IS}_k = \textbf{h} \left( \mathbf{f}_m^k(x_i) \right), \textrm{ where } \mathbf{f}_m^k(x_i) = \{f(x_j): x_j \in S_{m}^k(x_i) \},
\end{equation}
where $\textbf{h}: \mathbb{R}^{2m+1} \rightarrow \mathbb{R}$ is a polynomial of order 2 in stencil values. The order 2 is the direct consequence from using the $L^2$ norm in the definition for smoothness indicator.
\end{prop}

\begin{coro} [\textbf{The weights $\bm{\omega}_k$ are rational functions in the stencil values}] \label{coroA1} For a positive real number $\varepsilon > 0$, positive integer power $p$, and $\mathcal{C}_k^r \in \mathbb{R}^{+}$, the weights $\bm{\omega}_k$ in Eq.~\ref{complex_weights} are rational functions in the stencil values $\mathbf{f}_m(x_i)$ following proposition \ref{discrete_SI}, i.e.,
\begin{align}
    \bm{\omega}_k &= \frac{\mathcal{C}_k^r \prod_{l=0, l\neq k}^{r-1} (\varepsilon + \mathbf{IS}_l)^p }{\sum_{l=0}^{r-1} C_l^r \prod_{q=0,q\neq l}^{r-1}(\varepsilon + \mathbf{IS}_q)^p } \\
    &= \frac{\textbf{g}_1 (\mathbf{f}(x_i))}{\textbf{g}_2 (\mathbf{f}(x_i))}, \textrm{ where } \mathbf{f}(x_i) = \{f(x_j): x_j \in \bigcup\limits_{k=0}^{r-1} S_m^k \}.
\end{align}
Here, $\textbf{g}_1: \mathbb{R}^{2m+1} \rightarrow \mathbb{R}$, $\textbf{g}_2: \mathbb{R}^{2m+1} \rightarrow \mathbb{R}$ are both multi-variate polynomials of maximum degree $2p(r-1)$. The denominator polynomial $\textbf{g}_2$ will always be strictly positive due to the numerical corrections from $\varepsilon > 0$.
\end{coro}

\begin{coro} [\textbf{WENO approximation with weights $\mathbf{\omega}_k$ is rational function in the stencil values}] For linear map $g(\cdot)$ in Eq.(\ref{linear_map}), the WENO reconstruction of the function $f$ at intermediate point $x_{i+\frac{1}{2}}$ is also a rational function in stencil values following corollary \ref{coroA1}.
\begin{equation} 
    \hat{f}_{i+\frac{1}{2}} = \sum_{k=0}^{r-1} \bm{\omega}_k g\left( f_{i+k-r+1},...,f_{i+k} \right) = \frac{\textbf{p}_1 (\mathbf{f}(x_i))}{\textbf{p}_2 (\mathbf{f}(x_i))}, \textrm{ where } \mathbf{f}(x_i) = \{f(x_j): x_j \in \bigcup\limits_{k=0}^{r-1} S_m^k \}.
\end{equation}
Here, $\textbf{p}_1: \mathbb{R}^{2m+1} \rightarrow \mathbb{R}$, $\textbf{p}_2: \mathbb{R}^{2m+1} \rightarrow \mathbb{R}$ are multi-variate polynomials of maximum degree $2pr-1$ and $2p(r-1)$, respectively.
\end{coro}

\textbf{NOTE:} For the ease of analysis and presentation, we do not account for flux-splitting that helps with automatic up-winding in the above formulations. However, we note that the above analysis holds true for any smooth global flux splitting methods.

\section{Benchmark problems and training data generation} \label{data_gen_append}

We present the mathematical models used for the test cases and detail the procedure of data generation.

\subsection{Forced Burgers' case}
The forced Burgers' equation in 1D is given by the PDE:
\begin{equation}
    \frac{\partial u}{\partial t} + \frac{\partial (u^2)}{\partial x} = D \frac{\partial^2 u}{\partial x^2} + f(x,t)
\end{equation}
for the unknown function $u(x,t)$ with diffusion constant $D=0.02$. Here, we use the forcing term
\begin{equation}\label{forcing_parameters_burgers}
    f(x,t) = \sum_{i=1}^{N} A_i \sin(\omega_i t + 2\pi l_i x/L + \phi_i)
\end{equation}
with each parameter drawn independently and uniformly at random from its respective range: $A = [ -0.1,0.1 ]$, $\omega = [-0.4,0.4]$, $\phi = [0, 2\pi]$, and $N=20$. The domain size $L$ is set to $2\pi$ (i.e., $x\in [0,2\pi]$) with periodic boundary conditions, and $l = \{2,3,4,5\}$. We use a smooth initial condition $u(x,0) = \exp({-(x-3)^2}) $.

We generate the simulation data using finite differences on $N_x = 256$ evenly spaced grid points with fifth-order WENO discretization of the convection term. Time integration is performed using a third-order Runge-Kutta method with time-step size $\Delta t$  chosen as large as possible according to the CFL condition.

For larger domains we adjust the range of $l$ to preserve the wavenumber range. E.g., for $L = 8\pi$ we use $l = \{8,9,\ldots ,40\}$. We use the same spatial resolution $\Delta x$ for all domain sizes, i.e., the total number of grid points $N_x$ is adjusted proportionally to domain size.

\textbf{Training data conditions} \\
For training we use a coarse-graining factor $C = \{2,4,8\}$ in space. We sub-sample the data by removing  intermediate points in order to achieve the correct  coarse-graining factor $C$. For example, for $C=2$, we only keep spatial grid points with even indices. During training, time integration is done with a step size that satisfies the CFL-condition, i.e., $\Delta t_c \leq (\Delta x_c)^2/D$, where $\Delta x_c = C\Delta x$ and $\Delta t_c$ are the coarse spatial and time resolutions, respectively. The training domain in time is chosen until time $40$.

\subsection{Kuramoto-Sivashinsky (KS) case}
The Kuramoto-Sivashinsky (KS) equation for an unknown function $u(x,t)$ in 1D is given by the PDE:
\begin{equation}
    \frac{\partial u}{\partial t} + \frac{\partial (u^2)}{\partial x} + \frac{\partial^2 u}{\partial x^2} +\frac{\partial^4 u}{\partial x^4} = 0.
\end{equation}
The domain $x\in [-32,32]$ of length $L=64$ is discretized by $N_x = 256$ evenly spaced spatial grid points.

We use the initial condition
\begin{equation}
    u(x,t=0) = \sum_{i=1}^{N} A_i \sin(2\pi l_i x/L + \phi_i).
\end{equation}
Each parameter is drawn independently and uniformly at random from its respective range: $A = [ -0.5,0.5 ]$, $\phi = [0,2\pi]$ for $l = \{1,2,3\}$. 
We use a spectral method to numerically solve the KS equation using the \textsc{chebfun} (\cite{driscoll2014chebfun}) package with periodic boundary conditions. Time integration is performed using a fourth-order stiff time-stepping scheme with step-size $\Delta t = 0.05$.

\textbf{Training data conditions} \\
For the results presented in section \ref{KS_section}, we use a coarse-graining factor of $C=4$ in space, i.e., $N_c=64$, and we train STENCIL-NET on data up to $50$ time.
The time-step size is chosen as large as possible to satisfy the CFL condition.

\subsection{Korteweg-de Vries (KdV) case}

The Korteweg-de Vries (KdV) equation for an unknown function $u(x,t)$ is given by the PDE
\begin{equation}
    \frac{\partial u}{\partial t} + \frac{\partial (u^2)}{\partial x} + \delta \frac{\partial^3 u}{\partial x^3}  = 0,
\end{equation}
where we use $\delta = 0.0025$. We again use a spectral method to numerically solve the KdV equation using the \textsc{chebfun} (\cite{driscoll2014chebfun}) package. The domain is $x \in [ -1,1 ]$ with periodic boundary conditions, and the initial condition is $u(x,0) = \cos(\pi x)$. 

The spectral solution is represented on $N_x = 256$ equally spaced grid points discretizing the spatial domain.

\textbf{Noise addition details} \\
We corrupt a solution vector $\mathbf{U} \in \mathbb{R}^{N}$ with additive Gaussian noise as:
\begin{equation}
    \mathbf{V} = \mathbf{U} + \eta ,
\end{equation}
with $\eta = \sigma \cdot \mathcal{N}(0, \textrm{std}^2(\mathbf{U}))$. $\mathcal{N}(m,V)$ is the standard normal distribution with mean $m$ and variance $V$, and $\textrm{std}(\mathbf{U})$ is the empirical standard deviation of the solution $\mathbf{U}$. The parameter $\sigma$ is the magnitude of the noise.

\textbf{Training data conditions} \\
For the de-noising experiment presented in section \ref{KDv_section}, the value of $\sigma$ is set to 0.3. The spectral solution is sub-sampled with $C = 8$ in space and the time-step size $\Delta t = 0.02$ during training. We use a third-order TVD Runge-Kutta method for time integration during optimization.

\section{Optimization formulation}
\subsection{Space and time discretization of PDEs}\label{PDE_optim_appendix}
The continuous PDE in Eq.~\ref{PDE} is discretized in space and time. Following the discretization methods presented in Eqs.~\ref{linear_derv} and \ref{WENO_form}, we compactly formulate the spatially discretized version of the PDE as the system of ODEs
\begin{equation}\label{discrete_space_PDE}
    \frac{\partial u_i}{\partial t} = \mathcal{N}_d \big( \mathbf{u}_m(x_i),\, \Xi, \Delta x \big) + \mathcal{O}(\Delta x^{r_1}), \quad i=1,\ldots ,N_x,
\end{equation}
where $u_i = u(x_i)$, $\mathbf{u}_m(x_i) = \{u(x_j): x_j \in S_m(x_i) \}$, $\Delta x$ the grid resolution of the spatial discretization, and $\mathcal{N}_d: \mathbb{R}^{2m+1} \rightarrow \mathbb{R}$ is the nonlinear discrete approximation of the continuous nonlinear right-hand side operator $\mathcal{N}(\cdot)$.  The discrete approximation converges to the continuous operator with {\em spatial convergence rate} $r_1$ as $\Delta x \to 0$. Popular examples of spatial discretization methods include finite-difference, finite-volume, and finite-element methods. For simplicity of illustration, we only consider the case $\alpha_1 = 1$, $\alpha_2 =0$; the other case is analogous.

Neglecting the approximation error and integrating Eq.~\ref{discrete_space_PDE} on both sides over one time step $\Delta t$ yields a discrete map from $u(x_i,t)$ to $u(x_i,t+\Delta t)$ as:
\begin{equation}
    u_i^{n+1} = u_i^{n} + \displaystyle\int\limits_{t}^{t+\Delta t}   \mathcal{N}_d \big( \mathbf{u}_m^{n}(x_i), \Xi,  \Delta x \big) \, d\tau 
\end{equation}
where $u_{i}^{n+1} = u(x_i,t+\Delta t)$, $u_i^n = u(x_i, t)$, and $\mathbf{u}_m^n (x_i) = \{u(x_j,t): x_j \in S_m (x_i) \}$. Approximating the integral on the right-hand side of this map by quadrature, we find
\begin{equation}
    u_i^{n+1} = \mathbf{T}_d \bigg( u_{i}^n, \, \mathcal{N}_d \big( \mathbf{u}_m^n(x_i), \Xi, \Delta x \big),  \Delta t\bigg) + \mathcal{O}(\Delta t^{r_2}), \quad n=0,\ldots ,N_t -1,
\end{equation}
where $N_t$ is the total number of time steps.
Here, $\mathbf{T}_d $ is the explicit discrete time integrator with time-step size $\Delta t$. Due to approximation of the integral by quadrature, the discrete time integrator converges to the continuous-time map with {\em temporal convergence rate} $r_2$ as $\Delta t \to 0$.  Popular examples of explicit time-integration schemes include forward Euler, Runge-Kutta, and TVD Runge-Kutta methods. In this work, we only consider Runge-Kutta-type methods and their TVD variants as described by~\cite{shu1998essentially}. 

Based on Corollary~\ref{stencilnet_coro}, we can effectively approximate the discrete nonlinear function $\mathcal{N}_d (\cdot)$ with MLPs, which leads to:
\begin{equation} \label{approx_theta}
    \hat{u}_i^{n+1} = \mathbf{T}_d \bigg( \hat{u}_{i}^n,  \mathcal{N}_{\theta}\big( \hat{\mathbf{u}}_m^n(x_i), \Xi \big),  \Delta t\bigg) ,
\end{equation}
where $\mathcal{N}_{\theta}: \mathbb{R}^{2m+1} \rightarrow \mathbb{R}$ are the local nonlinear mlpconv layer rules parameterized with weights $\theta$. 

On the assumption of point-wise corruption of the solution $u$ with i.i.d.~noise $\eta$, i.e., $v_i^n = u_i^n + \eta_i^n$, the above Eq.~\ref{approx_theta} can be extended for mapping the measurements from $v_i^{n}$ to $v_i^{n+k}$, as:
\begin{equation} \label{k_int}
    \hat{v}_i^{n+k} = \mathbf{T}_d^{k} \Big( v_i^n - \eta_i^n,  \mathcal{N}_{\theta}\big( \mathbf{v}_m^n(x_i) - \hat{\mathbf{n}}_m^n(x_i) , \Xi \big),  \Delta t\Big) + \hat{\eta}_i^{n+k},
\end{equation}
where $\mathbf{v}_m^n(x_i) = \{v(x_j,t): x_j \in S_m (x_i) \}$ are the noisy measurement data and $\hat{\mathbf{n}}_m^n(x_i) = \{\hat{\eta}(x_j,t): x_j \in S_m (x_i) \}$ the noise estimates on the stencil $S_m$ centered at point $x_i$ at time $t$. The superscript $k$ in $\mathbf{T}_d^k$ denotes that the map integrates $k$ time steps into the future.

\subsection{Loss function and training} \label{optim}
From the Eq.~\ref{k_int}, we formulate a loss function for learning the local nonlinear discretization rules $\mathcal{N}_{\theta}$:
\begin{equation}\label{MSE}
      \mathcal{L}_{MSE} =  \sum_{n=1}^{N_t} \sum_{i=1}^{N_x} \sum_{k=-q}^{q} \gamma_k \Big\Vert v_{i}^{n+k} - \bigg( \mathbf{T}_d^{k} \Big( v_i^n - \eta_i^n,  \mathcal{N}_{\theta}\big( \mathbf{v}_m^n - \hat{\mathbf{n}}_m^n , \Xi \big),  \Delta t\Big) + \hat{\eta}_i^{n+k} \bigg) \Big\Vert_{2}^{2} 
\end{equation}
The point-wise noise estimates $\hat{\eta}_i$ are modeled as latent quantities, which is the reason why STENCIL-NET can decompose dynamics from noise. The scalar $\gamma_k$ are decaying weights that account for the accumulating prediction error (\cite{rudy2019deep}). The positive integer $q$ is the number of Runge-Kutta integration steps considered during optimization, which we refer to as the {\em training time horizon}. We penalize the noise estimates $\hat{\mathbf{N}}$
in order to prevent learning the trivial solution associated with the minimization problem given by Eq.~(\ref{MSE}). We also impose penalty on the weights of the network $\mathbf{W}$ in order to prevent over-fitting. The total loss is then given by
\begin{equation}
    \textrm{Loss} = \mathcal{L}_{MSE} + \lambda_n \Vert \hat{\textbf{N}}\Vert_F^2 + \lambda_{wd} \sum_{i=1}^{l} \Vert \mathbf{W}_i \Vert_F^2,
\end{equation}
where $l$ is the number of network layers, $\hat{\mathbf{N}} \in \mathbb{R}^{N_x \times N_t}$ is the matrix of point-wise noise estimates, and $\Vert \cdot \Vert _F$ is the Frobenius norm of a matrix. 

We use grid search through hyper-parameter space in order to identify values of the penalty parameters. We find the choice  $\lambda_n = 10^{-5}$ and $\lambda_{wd} = 10^{-8}$ to work well for all problems considered in this work. Alternatively, methods like drop-out, early-stopping criteria, and Jacobi regularization can be used to counter the over-fitting problem. Training is done using an Adam optimizer (\cite{kingma2014adam}) with learning rate $lr = 0.001$. For the activation we prefer the Exponential Linear Unit (ELU) for its smooth function interpolation abilities. One could readily use ReLU, Leaky-ReLU, \textrm{tanh} or adaptive activation functions with learnable hyper-parameters, however, we choose ELU for its superior empirical performance for our class of problems.

\textbf{code availability:} We provide STENCIL-NET implementation in PyTorch with in-built finite-difference solvers and WENO module for reproducing the simulation data from forced Burgers' problem. We also include the  spectral solution used for training in the KS and KdV problems. We provide jupyter notebook demos with pre-trained models for reproducing some results from the main text. To preserve the anonymity of the submission we provide our codes as additional files in the supplementary material.

\newpage

\section{Additional results on accuracy and the choice of hyper-parameters}  \label{more_results_appendix}
\setcounter{figure}{0} \renewcommand{\thefigure}{D.\arabic{figure}}
We present detailed results on the accuracy of STENCIL-NET. We use the mean square error (MSE) metric to present the convergence of STENCIL-NET solutions:
\begin{equation*}
    \textrm{MSE} =   \frac{1}{N_t N_x}\sum_{n=1}^{N_t} \sum_{i=1}^{N_x} \big \Vert \hat{u}_i^n - u(x_i,t)\big \Vert_{2}^{2}.
\end{equation*}
Here, $\hat{u}$ is the STENCIL-NET prediction and $u(x_i,t)$ the true solution value at grid point $x_i$ and time $t=n\Delta t$. 

\begin{figure}[!htbp]
    \centering
    \includegraphics[width=1.0\textwidth]{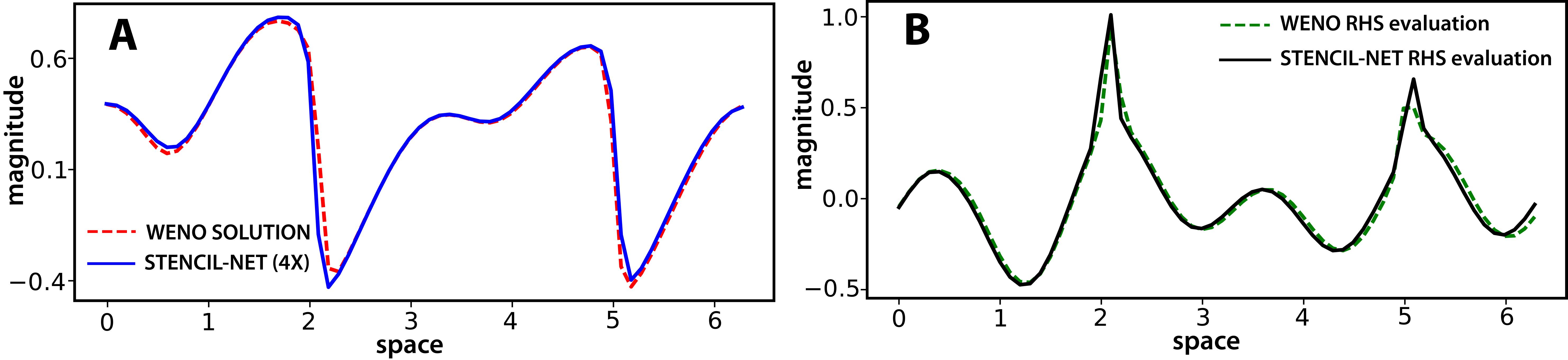}
    \caption{\small\textbf{Comparison between STENCIL-NET (4-fold coarsened) and fifth-order WENO for forced Burgers'}: ({A}) Comparison between the STENCIL-NET prediction and WENO solution over the entire domain at time $t= 40$. ({B}) Comparison between the corresponding nonlinear discrete operators (PDE RHS) computed by STENCIL-NET ($\mathcal{N}_{\theta}$) and by the WENO scheme ($\mathcal{N}_{d}$) at time $t=40$.}
    \label{burgers_time_slice}
\end{figure}

\begin{figure}[!htbp]
    \centering
    \includegraphics[width=1.0\textwidth]{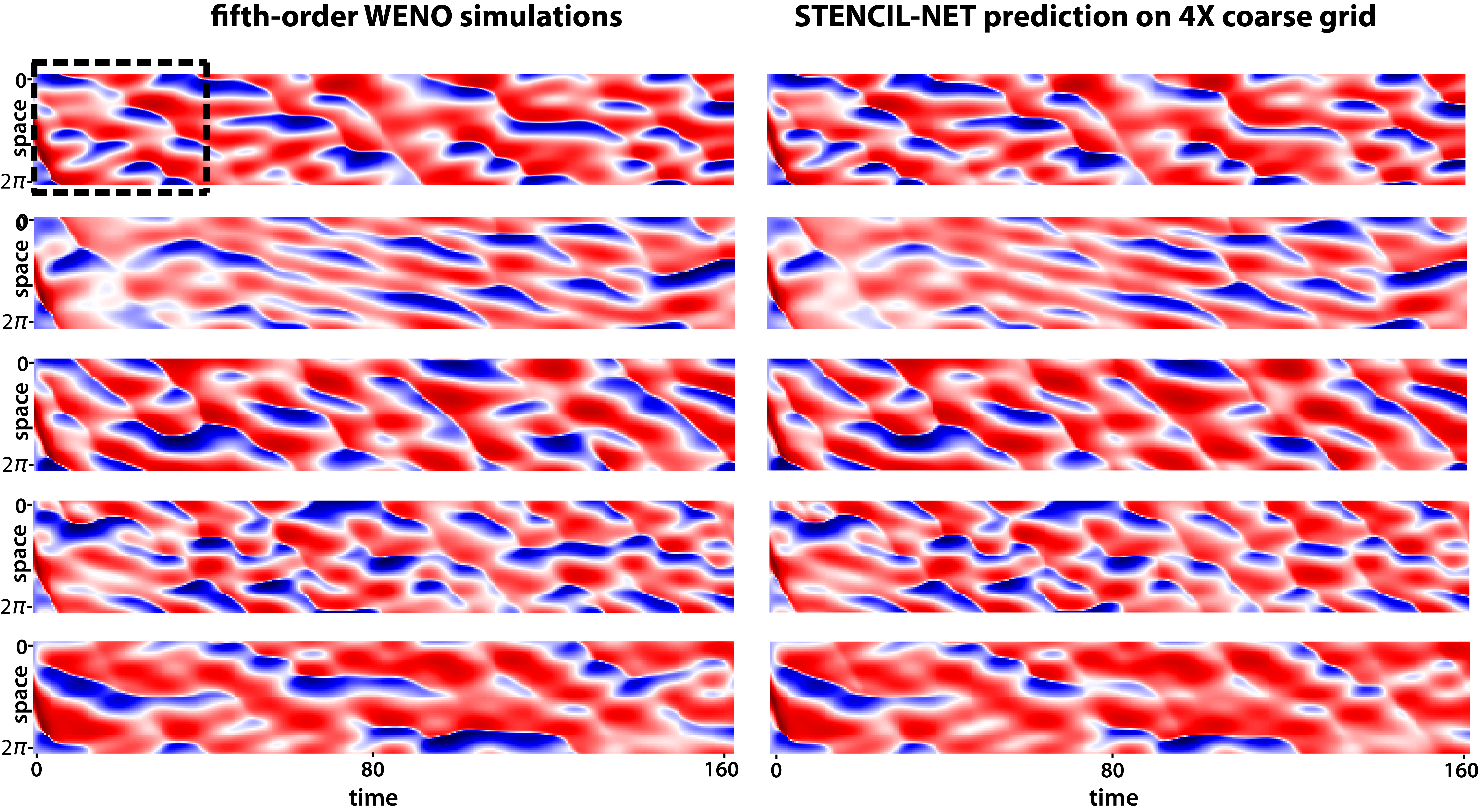}
    \caption{\small \textbf{STENCIL-NET prediction for different forcing terms of forced Burgers' without  re-training} Each row corresponds to a forced Burgers' solution with different forcing term parameters (see Eq.~\ref{forcing_parameters_burgers}). 
    STENCIL-NET was trained only once using data from the forcing term in the top row (training domain marked by dashed box).
    The left plots show the fifth-order WENO solutions for the different forcing terms, the right plots are the STENCIL-NET prediction on $4\times$ coarser grids without additional training. 
    }
    \label{general_forcing}
\end{figure}

STENCIL-NET relies on mapping the solution on local patches to rich discretization features. This allows STENCIL-NET to resolve  spatio-temporal features that contain steep gradients and discontinuities as shown in Fig.~\ref{burgers_time_slice}. The learned local nonlinear discretization rules can then be used to generalize for parameters and domain-sizes beyond the training conditions. This is illustrated in Fig.~\ref{general_forcing} where a STENCIL-NET model that has been trained on a small training domain can generalize well to longer times and to equations containing different forcing terms $f(x,t)$ in Eq.~\ref{forcing_parameters_burgers}.

\begin{figure}[!htbp]
    \centering
    \includegraphics[width=1.0\textwidth]{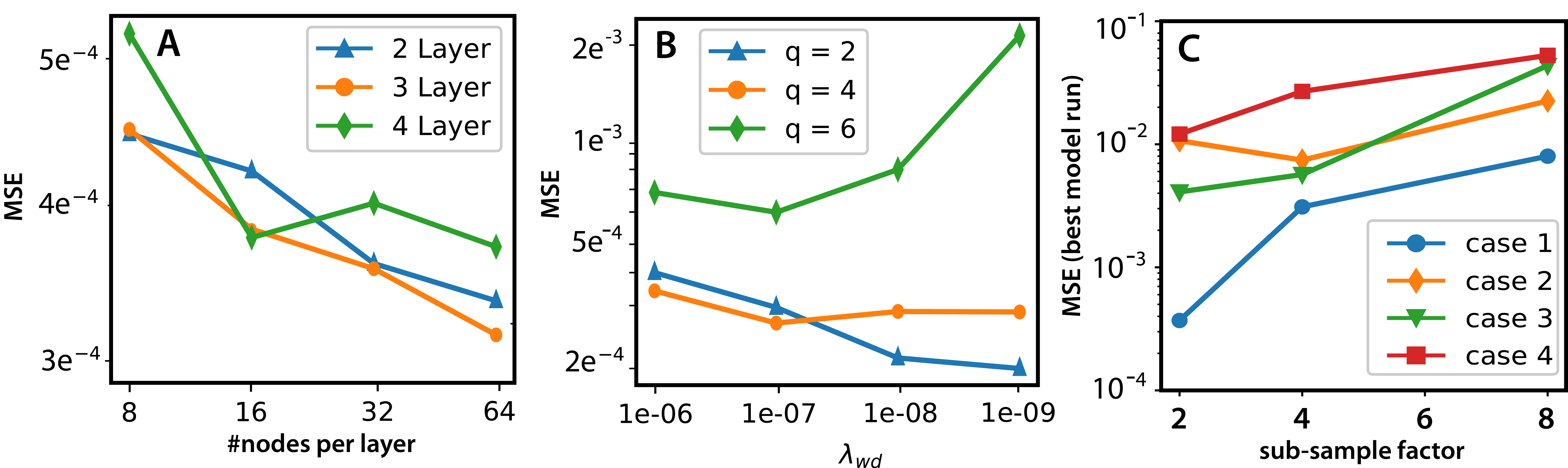}
    \caption{\small \textbf{Architecture choice, choice of $\boldsymbol{\lambda_{wd}}$, and generalization power}: All the models are trained on the spatio-temporal data contained in the dashed box in Fig.~\ref{burgers_prediction} of the main text. The training box encompasses the entire spatial domain of length $L=2\pi$ and a time duration of $T=40$.  ({A)} STENCIL-NET prediction accuracy for different choices of MLP architectures. Each data point corresponds to the average MSE accumulated over all the stable seed configurations.
    ({B)} Effect of the weights regularization parameter $\lambda_{wd}$ on the prediction accuracy for different numbers of Runge-Kutta integration steps $q$ during optimization. The plots in B were produced with a 3-layer, 64-nodes MLP architecture, which showed the best performance in A.
    ({C)} Generalization power to larger domains: Case 1: STENCIL-NET MSE for predictions on space and time domains of the same extent as the training box. Case 2: STENCIL-NET prediction MSE on the same domain length ($L=2\pi$) but for $4\times$ longer time duration than the training box, i.e., $L=2\pi$, $T=160$. Case 3: STENCIL-NET prediction MSE for $16\times$ larger spatial domain but the same time duration as the training box, i.e., $L=32\pi$, $T=40$. Case 3: STENCIL-NET prediction MSE for $16\times$ larger spatial domain and $4\times$ longer time duration than the training box, i.e., $L=32\pi$, $T=160$. We use the best model specific to the resolution (sub-sample factor) with the smallest prediction error to produce each data point in the plot. } 
    \label{hyperparams}
\end{figure}

In Fig.~\ref{hyperparams}A,B, we show the effect of the choice of MLP architectures and of the weight regularization parameter $\lambda_{wd}$ for different training time horizons $q$. We  quantify the generalization power of STENCIL-NET in Fig.~\ref{hyperparams}C for predictions of varying resolutions and for different lengths of the space and time domains.

\begin{figure}[h]
    \begin{minipage}[b]{1.0\textwidth}
    \centering
    \includegraphics[width=0.6\textwidth]{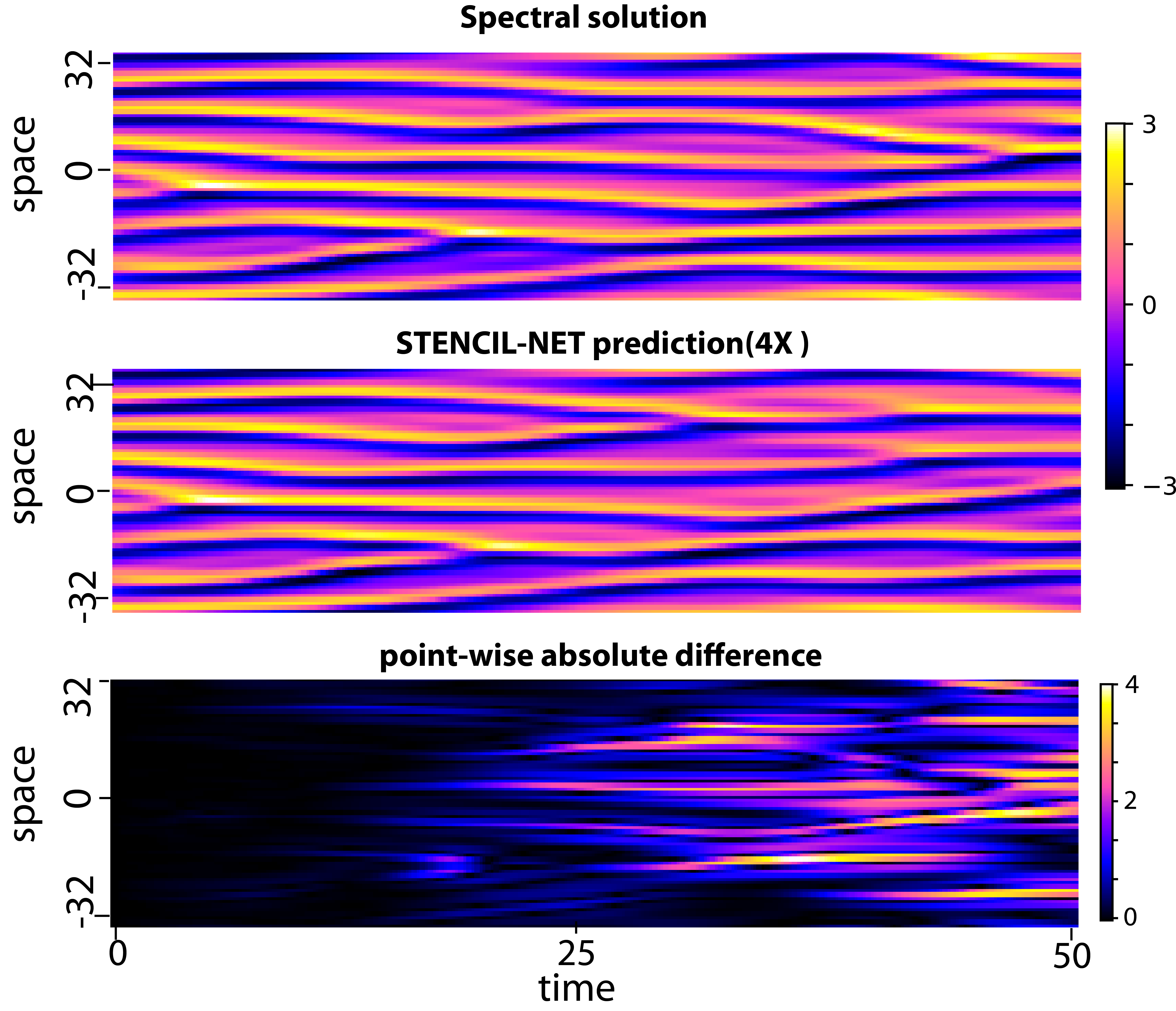}
    \end{minipage}
    \caption{\small \textbf{Equation-free prediction of chaotic spatio-temporal KS dynamics:} (Top) Spectral solution of the KS equation on domain $L=64$. (Middle) STENCIL-NET prediction on a $4\times$ coarser grid. (Bottom)  
    Point-wise absolute error between the spectral and STENCIL-NET solutions.}
    \label{KS_prediction}
\end{figure}

The spatio-temporal dynamics of the chaotic KS system is shown in Fig.~\ref{KS_prediction}. In contrary to the forced Burgers' case, the STENCIL-NET predictions on a $4\times$ coarser grid start to diverge from the true spectral solution (see bottom row of Fig.~\ref{KS_prediction}). This is due to the chaotic behavior of the KS system for space domain lengths $L>22$, which causes any small prediction error to grow exponentially at a rate proportional to the maximum Lyaponuv exponent of the system. Despite this fundamental unpredictability of the actual solution, is able to correctly estimate the value of the maximum Lyaponuv exponent and the spectral statistics of the KS system (see Fig.~\ref{KS_metrics} in the main text). This is evidence that the equation-free STENCIL-NET model has correctly learned  the intrinsic ergodic properties of the system that has generated the data.

\begin{figure}[h]
\begin{minipage}[b]{1.0\textwidth}
\centering
\includegraphics[width=1.0\textwidth]{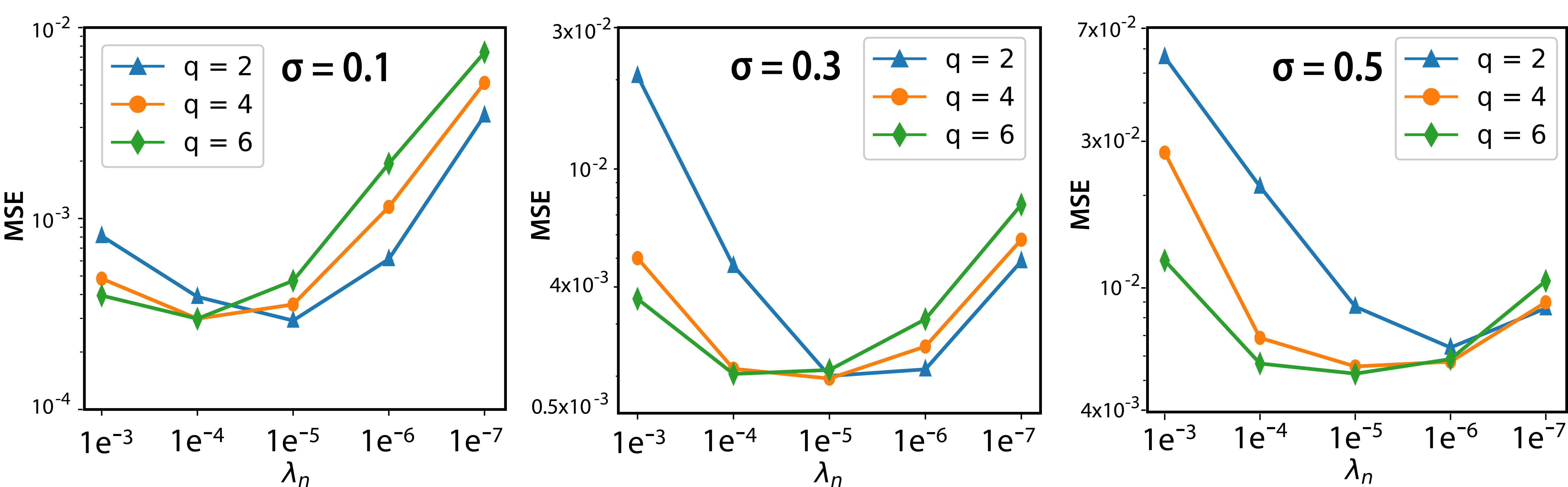}
\end{minipage}
    \caption{\small \textbf{Effect of the noise-penalty $\boldsymbol{\lambda_{n}}$ on de-noising:} We show the STENCIL-NET prediction MSE for different noise penalty $\lambda_n$ and training time horizon $q$. The data was corrupted by additive Gaussian noise with $\sigma=0.1,\, 0.3, \, 0.5$ (from left to right).}
    \label{KDV_noise}
\end{figure}

For the de-noising case, we report in Fig.~\ref{KDV_noise} the effect of the noise-penalty parameter $\lambda_n$ for varying noise levels. We also report for each case the dependence on the training time horizon $q$. This empirical evidence suggests that the choice of $q$ can be used to trade off between stability and accuracy of the prediction.

\newpage

\section{Estimation of the reduction in computational elements and the computational speed-up} \label{speedup_appendix}

We derive an estimate for the speed-up achieved by replacing a numerical solver with STENCIL-NET on a $C$-times coarser grid. For the sake of argument, we only consider problems that contain a diffusion term in the PDE model. For those, the time-step size $\Delta t$ of an explicit time integrator is bounded by the CFL condition.
\begin{equation*}
    \Delta t \leq \frac{(\Delta x)^2}{2D} .
\end{equation*}
Such a CFL condition also exists for other, non-diffusive problems, but it may take a different form. The estimation below should then be adjusted accordingly. 

In any case, the number of computational grid points per dimension is inversely proportional to the resolution, i.e. $N_x \propto 1/\Delta x$. This implies that the total number of computational elements (i.e., grid points or particles) or degrees of freedom (DOF) needed by a numerical scheme in $d$ dimensions is
\begin{equation*}
\textrm{DOF} = \frac{K}{(\Delta x)^d \cdot \Delta t},
\end{equation*}
where $K$ is a positive constant that depends on the domain size $L$ and the total time $T$ the simulation is run for. Coarse-graining the spatial-domain by a factor of $C>1$ leads to $\Delta x_c = C \Delta x$. According to the CFL condition, this allows taking larger time steps $\Delta t_c = C^2 \Delta t$ and therefore reduces the number of degrees of freedom to:
\begin{equation}\label{speedup_factor}
    \textrm{DOF}_c = \frac{K}{ (C\Delta x)^d \cdot  (C^2\Delta t )  } = \frac{\textrm{DOF}}{C^{d+2}}.
\end{equation}
Thus, coarse-graining using methods like STENCIL-NET can lead to a drastic reduction in computational elements, with little compromise on accuracy and cost per time-step. For example, being able to coarse-grain by a factor of $C=8$ in a 3D simulation would reduce the DOF 32,768-fold. These analysis also hold true for the Euler systems which have no diffusion term, however, the number of degrees of freedom now become,
\begin{equation*}
    \textrm{DOF}_c = \frac{K}{ (C\Delta x)^d \cdot  (C^2\Delta t )  } = \frac{\textrm{DOF}}{C^{d+1}}.
\end{equation*}
This allow us to readily employ STENCIL-NET for speeding the solution of Euler like systems and level-set formulations.

Let $t_n$ and $t_s$ be the computational time required for each DOF in a numerical method and in STENCIL-NET, respectively. Then from Eq.(\ref{speedup_factor}), the solution-time speed-up achieved by STENCIL-NET substitution is
\begin{equation*}
    \kappa = \left( \frac{t_s}{t_n} \right) C^{d+2}.
\end{equation*}
The ratio $(t_n/t_s)$ is the {\em computational overhead} of STENCIL-NET. As shown in Fig.~\ref{speedup} of the main text, this ratio is between 1 and 2 on GPUs, and between 14 and 20 on CPUs. 
This, for example, leads to speed-ups $\kappa = 16384\ldots 32768$ on the GPU or $\kappa = 1638\ldots 2340$ on the CPU for an 8-fold coarse-grained simulation in 3D. Given that real-time applications, for example in virtual reality, require frame rates of about 60\,Hz, this implies that even a sequential CPU implementation of STENCIL-NET should be able to achieve real-time performance for PDE simulations that would normally require about 30 seconds per time step of computer time. This is well within the range of practical applications.

\end{document}